\documentclass[11pt,reqno]{amsart}
\usepackage{amsthm,amssymb,amsbsy,epsfig,graphicx,color}
\usepackage{hyperref}

\input{amssym.def}
\input{amssym}

\def\be{\begin{equation}}
\def\ee{\end{equation}}
\def\ba{\begin{eqnarray}}
\def\ea{\end{eqnarray}}
\def\lb{\label}
\def\nn{\nonumber}
\def\ni{\noindent}

\def\cal{\mathcal}

\def\R{{R}}
\def\Rd{{\stackrel{\ast}{R}}}
\def\Jd{{\raisebox{1pt}{$\stackrel{\raisebox{-1pt}{$\scriptstyle\,\, \ast$}}{J}$}}}
\def\id{{I}}

\def\O{{ O_{\hspace{-2pt}\raisebox{-2pt}{\tiny\sl R}}}}
\def\D{{ D_{\hspace{-2pt}\raisebox{-2pt}{\tiny\sl R}}}}
\def\C{{ C_{\hspace{-2pt}\raisebox{-2pt}{\tiny\sl R}}}}
\def\Dm{{ D_{\hspace{-2pt}\raisebox{-2pt}{$\scriptscriptstyle R^{-1}$}}}}
\def\Cm{{ C_{\hspace{-2pt}\raisebox{-2pt}{$\scriptscriptstyle R^{\scriptscriptstyle -1}$}}}}
\def\Dd{{ D_{\hspace{-1pt}\stackrel{\raisebox{-1pt}{$\scriptscriptstyle\,\ast$}}
{\raisebox{1pt}{$\scriptscriptstyle R$}}}}}

\def\Cd{{ C_{\hspace{-1pt}\stackrel{\raisebox{-1pt}{$\scriptscriptstyle\,\ast$}}
{\raisebox{1pt}{$\scriptscriptstyle R$}}}}}

\def\Rdro{{ \rho_{\hspace{-1pt}\stackrel{\raisebox{-1pt}{$\scriptscriptstyle\,\ast$}}
{\raisebox{1pt}{$\scriptscriptstyle R$}}}}}

\def\Rro{{\rho_{\hspace{-1pt}\raisebox{-2pt}{\tiny\sl R}}}}
\def\Rpsi{{\Psi_{\hspace{-2pt}\raisebox{-2pt}{\tiny\sl R}}}}
\def\tr{{\rm Tr}\,}
\def\Adro#1{\stackrel{\raisebox{-1pt}{ $\scriptstyle \ast$}}
{A}^{\raisebox{-5pt}{$\scriptstyle #1$}}}

\def\Md#1{M_{\!\stackrel{\raisebox{-5pt}{$\scriptscriptstyle\hspace{.5pt}\ast$}}
{\raisebox{5pt}{$\scriptscriptstyle \overline{#1}$}}}}

\def\subRd{\!{\raisebox{1pt}{$\stackrel{\raisebox{-1pt}
{$\scriptscriptstyle\,\, \ast$}}{\scriptstyle R}$}}}

\def\Rtr{{\rm Tr}_{\hspace{-2pt} \raisebox{-2pt}{\tiny\sl R}\,}}
\def\RTr#1{{\rm Tr}_{\hspace{-2pt} \raisebox{-2pt}{\tiny\sl R}\,\raisebox{2pt}{\scriptsize$(#1)$}}}
\def\Tr#1{{\rm Tr}_{\raisebox{2pt}{\scriptsize$(#1)$}}}
\def\Rdtr{{\rm Tr}_{\hspace{-2pt} \raisebox{-2pt}{{\raisebox{1pt}
{$\stackrel{\raisebox{-1pt}{$\scriptscriptstyle\,\, \ast$}}{\scriptstyle R}$}}}\,}}

\def\RdTr#1{{\rm Tr}_{\hspace{-2pt} \raisebox{-2pt}{\raisebox{1pt}
{$\stackrel{\raisebox{-1pt}{$\scriptscriptstyle\,\, \ast$}}
{\scriptstyle R}$}}\,\raisebox{2pt}{\scriptsize$(#1)$}}}

\def\Rdet{{{\rm det}_{\hspace{-1pt} \raisebox{-2pt}{\tiny\sl R}\,}}}
\def\RX#1{X_{\hspace{-2pt} \raisebox{-2pt}{\tiny\sl R}\,}^{#1}}

\makeatletter\@addtoreset{equation}{subsection}\@addtoreset{equation}{section}\makeatother

\newcounter{theorem}
\makeatletter\@addtoreset{theorem}{section}\makeatother

\newtheorem{prop}[theorem]{Proposition}
\newtheorem{rem}[theorem]{{\em Remark}}
\newtheorem{lem}[theorem]{Lemma}
\newtheorem{def-lem}[theorem]{Definition-Lemma}
\newtheorem{def-prop}[theorem]{Definition-Proposition}
\newtheorem{defin}[theorem]{Definition}
\newtheorem{theor}[theorem]{Theorem}
\newtheorem{cor}[theorem]{Corollary}
\newtheorem{example}[theorem]{{\em Example}}

\begin{document}
\begin{titlepage}
\title[Spectral extension of the quantum group cotangent bundle]
{Spectral extension of the quantum group\\ cotangent bundle
%:\\\vspace{1.5mm}
%{\smaller\smaller calculating dynamical R-matrices and\\
%analyzing  evolution of the isotropic $q$-top}\vspace{0.0cm}
}

\author{Alexei Isaev}
\address{Alexei Isaev,
Bogoliubov Laboratory of Theoretical Physics, Joint Institute
for Nuclear Research, 141980 Dubna, Moscow Region, Russia}
\email{isaevap@theor.jinr.ru}

\author{Pavel Pyatov}
\address{Pavel Pyatov,
Max Planck Institute for Mathematics,
Vivatsgasse 7, D-53111 Bonn, Germany \&
Bogoliubov Laboratory of Theoretical Physics, Joint Institute
for Nuclear Research, 141980 Dubna, Moscow Region, Russia}
\email{pyatov@theor.jinr.ru}
%\thanks{}

\date{\today}

\begin{abstract}
The structure of a cotangent bundle is investigated for
quantum linear groups $GL_q(n)$ and $SL_q(n)$.
Using a $q$-version of the Cayley-Hamilton theorem
we construct an extension  of the algebra
of  differential operators on $SL_q(n)$ (otherwise called the Heisenberg double)
by  spectral values of the matrix of right invariant vector fields.
We consider two applications for the spectral extension.
First, we describe the extended Heisenberg double in terms of a new set of generators ---
the Weyl partners of the spectral variables. Calculating defining relations in terms of these
generators allows us to derive $SL_q(n)$ type dynamical R-matrices in a surprisingly simple way.
Second, we calculate an evolution operator for the model of
$q$-deformed isotropic top introduced by A.Alekseev and L.Faddeev.
The evolution operator is not uniquely defined and
we present two possible expressions for it. The first one is a Riemann theta function
in the spectral variables. The second one is an almost free motion evolution
operator in terms of logarithms of the spectral variables.
Relation between the two operators is given by a modular functional equation
for Riemann theta function.
\end{abstract}

\maketitle

\end{titlepage}

\phantom{a}
\vspace{2cm}
%\vfill
\tableofcontents
\vfill\newpage

\section{Introduction}

A notion of a Heisenberg double over quantum group has been formulated
and attracted substantial researcher's interest
in the early 90-s \cite{AF1,STS,SWZ1,SWZ2}. From the algebraic point of view it is a
smash product algebra (see \cite{Mon})
of the quantum group (or, the quantized universal enveloping algebra)
and its dual Hopf algebra (see \cite{Drin,FRT}).
In the differential geometric interpretation it may be viewed
as an algebra of quantized differential operators over group  or, equivalently, as
an algebra of quantized functions over cotangent bundle of the group.
Since the group's cotangent bundle serve a typical phase space for
integrable classical dynamics, it is natural to attach the same role to
the Heisenberg double over quantum group for quantum physical models.
As a test example, a model of   $q$-deformed isotropic
top was suggested in \cite{AF1,AF2}. A discrete time evolution in this model
is given by a series of automorphisms of the Heisenberg double. It turns out however
that finding an explicit expression for the model's
evolution operator is not just a technical problem.\footnote{
This problem was suggested to authors by L.D. Faddeev
in summer 1996, during Alushta conference ``Nonlocal, nonrenormalizable
field theories''.} The automorphisms defining the evolution by no means can be treated as
inner ones in the original algebra and so, for a proper realization of
the $q$-top one needs an appropriate extension of
the Heisenberg double.

\smallskip
Also stimulated by the invention of
quantum groups were general studies of the algebras whose generators
satisfy quadratic relations (see \cite{PP} and references therein) and investigations of minor
identities for matrices over noncommutative rings
\cite{GR1,GR2,KL}. These two lines of research are meeting together
in the theory of the so-called quantum matrix algebras \cite{H,IOP}
whose structure theory can be developed in a full analogy with the usual matrix
analysis. In particular, one can define quantum versions of the matrix trace and
determinant \cite{FRT}, introduce notions of a spectrum and a power of quantum matrix, and
formulate the Cayley-Hamilton theorem (see \cite{GPS,IOP,OP} and references therein).

\smallskip
A remarkable fact about quantum matrix algebras
is that their  most known examples --- the RTT algebra \cite{FRT} and the reflection equation
algebra \cite{KS} --- serve the building blocks for a construction of the
quantum group differential geometry in general \cite{SWZ1}
and so, also for the Heisenberg double.
It is the aim of the present paper to apply  the
structure results on the quantum matrix algebras
for investigation of the dynamics of the isotropic $q$-top.
Following \cite{AF1,STS,SWZ1} we begin with a definition of the Heisenberg double as a smash
product algebra of a pair of quantum matrix algebras. These are
the RTT algebra, playing the role of quantized functions over group,
and the reflection equation algebra, interpreted as
quantized right invariant differential operators over group.
We then consider a central extension of the reflection equation algebra by the
spectrum of it's generating matrix of quantized right invariant vector fields,
and define a proper (non central) extension
of the whole Heisenberg double by these spectral variables.
Finally, after the spectral extension is made, the evolution of the isotropic $q$-top becomes
an inner automorphism of the Heisenberg double.\footnote{Strictly speaking, one has to
extend the algebra by a formal power series in the spectral variables.} Constructing
the evolution operator is then straightforward.

\smallskip
The paper is organized as follows. In the next section we recall some facts about
universal R-matrix and the R-matrix techniques.
We are mainly discussing the case of  (numeric) R-matrices of a type $GL_q(n)$.
These type R-matrices are later on used for description of the cotangent bundles
(or, the Heisenberg doubles) over quantum linear groups.
% $GL_q(n)$ and $SL_q(n)$.

In section 3 we introduce the RTT algebra, the reflection equation algebra
and define their smash product algebra --- the Heisenberg double. We
are describing the algebras of the two linear types --- $GL_q(n)$ and $SL_q(n)$.
For the reflection equation algebra we formulate in these cases the Cayley-Hamilton
theorem and use it for the spectral extension of the Heisenberg double. This is
the first main result of the paper (see theorem \ref{main1}).

The Heisenberg double is initially defined in terms of the quantized right invariant vector fields.
In order to demonstrate the left-right symmetry of the Heisenberg double,
in subsection \ref{subsec3.4} we describe it using quantized
left invariant vector fields. We also derive explicit relations between the spectra
of the matrices of left and right  generators, see corollaries \ref{cor3.16} and \ref{cor3.33}.
To keep clearness of a presentation some technical lemmas are moved
from this subsection to appendix B.

The spectral extension suggests yet another distinguished generating set for the Heisenberg double, namely, the one
which satisfy the simplest possible -- Weyl algebraic -- relations with the spectral variables.
In the subsection \ref{subsec3.5} we derive defining relations for this set,
see theorem \ref{theor3.35}. Quite expectantly, the relations involve
dynamical R-matrices whose  dynamical arguments are
the spectral variables (see corollary \ref{cor3.36}). Surprising facts are that
the dynamical R-matrices are coming in pairs, and that they are derived by solving a
simple system of (at most three) linear equations.

%\smallskip
Section 4 is devoted to solving a dynamical problem for the isotropic $q$-top.
This is our second main result.
Noticing that an evolution operator of the model is not uniquely defined,
we derive two different expressions for it. The first one is given in terms of the
Riemann theta function whose matrix of periods is proportional to  Gram matrix of the
lattice $A^*_{n-1}$, see relations (\ref{add}), (\ref{i.14}). This solution converges
for $|q|<1$, or for $q$  a rational root of 1.
The second solution converging for arbitrary values of $q$ is given in terms of
logarithms of the spectral variables, see \eqref{i.x-0}, \eqref{i.x}.
The idea for the logarithmic substitution  (that means
passing from  Weyl type to Heisenberg type commutation relations)
was suggested to authors  by L.D. Faddeev (for argumentation see \cite{Fad,Faddd}).
The evolution in the logarithmic variables reduces to an almost free motion.
A relation between the two solutions is given then by a modular functional equation
for Riemann theta function  \eqref{Theta-1-2}.

\smallskip
Concluding the introduction we would like to mention a number of open problems which,
in our opinion, deserve further investigation. First of all, it is straightforward to
formulate a problem of spectral extension for the Heisenberg doubles over
orthogonal and symplectic quantum groups and over quantum linear supergroups.
Technical prerequisites for this were developed, respectively, in \cite{{OP}}
and \cite{GPS1,GPS2}. As well, it would be interesting to construct spectral extensions for the cases
of real quantum groups. We believe that a correct setting for this problem
is suggested in \cite{AF2}.

Another interesting problem is an extension of a modular double construction
\cite{Faddd,Fadub} (see also \cite{Leb1,Leb2}) for the case of Heisenberg double
over quantum group. A starting point for investigation here would be a modular functional
relation \eqref{Theta-1-2} between the two evolution operators constructed in section 4.
Riemann theta function standing
in the denominator in this relation could be considered as an evolution operator for the
modular dual Heisenberg double.

At last, an observation that  a ribbon element serves a $q$-top evolution operator
on the smash product algebra of a ribbon Hopf algebra with it's dual Hopf algebra
(see example \ref{example-ribbon})
could open a way for the spectral extension of a quasi-triangular Hopf algebra.
A partial step in this direction is made in appendix A, where pairing
%${\cal h}\cdot ,\cdot {\cal i}$
of the quasi-triangular Hopf algebra with its dual
Hopf algebra is extended for the set of spectral variables, see corollary \ref{cor-A2}.
\vspace{0.5cm}

{\bf Acknowledgement}.
We are grateful to Ludwig Dmitrievich Faddeev
for acquainting us with the problem of a dynamics of the isotropic q-top,
and for numerous inspiriting discussions and advises.
We would like to thank Alexei Gorodentsev, Sergei Kuleshov, Andrey Levin,
Dmitry Lebedev, Andrey Mudrov, Andrey Marshakov and  Vyacheslav Spiridonov
for their useful comments and conversations.
We also would like to acknowledge a warm hospitality of the Max-Planck-Institute
f\"{u}r Mathematik where writing this paper was started in 2004 and finished in 2008.
The work is supported by the
Russian Foundation for Basic Research,  grant No.\,08-01-00392-a.

\section{R-matrices}
\lb{sec2}

In this introductory
section we collect some necessary information about  R-matrices and
an R-matrix technique.

\subsection{Universal R-matrix}
\lb{subsec2.0}

First, we recall few basic notions  from the theory of
quasi-triangular  Hopf algebras \cite{Drin,Drinf} and ribbon Hopf algebras \cite{ReshT}
(for review see \cite{ChP,KSch}).

Let ${\goth A}$ be a Hopf $\mathbb C$-algebra supplied with a unit $1: {\mathbb C}\rightarrow {\goth A}$,
a counit $\epsilon: {\goth A}\rightarrow {\mathbb C}$,
a product $\textbf{m}:  {\goth A}\otimes {\goth A}\rightarrow {\goth A}$, a coproduct
$\Delta:{\goth A}\rightarrow {\goth A}\otimes {\goth A}$, and  an antipode
$S: {\goth A}\rightarrow {\goth A}$
mappings subject to standard axioms.

A Hopf algebra ${\goth A}$ is called {\em almost cocommutative} if there exists an invertible
element ${\cal R}\in {\goth A}\otimes {\goth A}$
that  intertwines the coproduct $\Delta$ and the opposite coproduct $\Delta^{op}$
(in Sweedler's notation:
$\Delta^{op}(x) = x_{(2)}\otimes x_{(1)}$ if $\Delta(x) = x_{(1)}\otimes x_{(2)}$)
\be
\lb{intertwine}
{\cal R} \Delta (x) = \Delta^{op}(x) {\cal R}\qquad \forall\, x\in{\goth A}_{\cal R}.
\ee
In this case the element ${\cal R}$ is called a {\em universal R-matrix}, and the corresponding
almost cocommutative  Hopf algebra is denoted as
${\goth A}_{\cal R}$.
The algebra ${\goth A}_{\cal R}$ is called {\em quasi-triangular}
if additionally $\cal R$ satisfies relations
\be
\lb{coco}
(\Delta\otimes id){\cal R} \,=\, {\cal R}_{13}{\cal R}_{23}\, , \quad
(id\otimes\Delta){\cal R} \,=\,  {\cal R}_{13}{\cal R}_{12}\, ,
\ee
where ${\cal R}_{12} = {\cal R}\otimes 1$,
${\cal R}_{23} = 1 \otimes {\cal R}$,  and ${\cal R}_{13}= \sum_i a_i\otimes 1\otimes b_i$
for ${\cal R}=\sum_i a_i\otimes b_i$.
Relations (\ref{intertwine}), (\ref{coco}) together imply an equality
%equalities
%\be
%\nn
%(\epsilon \otimes id)({\cal R})\, =\, (id\otimes\epsilon)({\cal R})\, =\,  1\, ,\qquad
%(S\otimes id)({\cal R})\, =\, (id\otimes S^{-1})({\cal R})\, =\,  {\cal R}^{-1}\, ,
%\ee
\be
\lb{UYB}
{\cal R}_{12} {\cal R}_{13} {\cal R}_{23} = {\cal R}_{23} {\cal R}_{13} {\cal R}_{12}\, ,
\ee
which is called the {\em Yang-Baxter equation}.

In the almost cocommutative case an element
$ u\, :=\, \textbf{m} (S\otimes id)({\cal R}_{21}) \in {\goth A}_{\cal R}$
is invertible. In terms of $u$ the square of antipode is expressed as
\be
\lb{square-antipode}
S^2(x)\, =\, u\, x\, u^{-1}\qquad \forall\, x\in {\goth A}_{\cal R}\, .
\ee
In the quasi-triangular case one has following  formulas
\ba
\nn
S(u)\, =\, \textbf{m} (id\otimes S)({\cal R}_{12})\, , &&
u^{-1}\, =\, \textbf{m} (id\otimes S^2)({\cal R}_{21})\, ,
\\[1mm]
\lb{Ruu}
\Delta(u)\, =\, ({\cal R}_{21}{\cal R}_{12})^{-1}\, u\otimes u\,
 , && {\cal R}\, (u\otimes u)\, =\, (u\otimes u)\, {\cal R}\, .
\ea

An element $ u S(u) = S(u) u$ (in the almost cocommutative case)
belongs to the center of ${\goth A}_{\cal R}$.
A {\em central} extension of the quasi-triangular Hopf algebra ${\goth A}_{\cal R}$
by a so-called {\em ribbon element} $\upsilon$ such that
\be
\lb{ribbon-def}
\upsilon^2\ =\, u S(u)\, , \qquad
\Delta (\upsilon)\, =\, ({\cal R}_{21}{\cal R}_{12})^{-1}\,\upsilon\otimes \upsilon
\ee
is called a {\em ribbon Hopf algebra}. The ribbon element also fulfills relations
\be
\nn
\epsilon(\upsilon)\, =\, 1\, , \qquad
S(\upsilon)\, =\, \upsilon\, , \qquad
{\cal R}\, (\upsilon\otimes\upsilon)\, =\,(\upsilon\otimes\upsilon)\,{\cal R}\, .
\ee

\smallskip
Throughout this paper our basic reference example of the quasi-triangular Hopf algebra
${\goth A}_{\cal R}$
is the  {\em quantized universal enveloping  algebra}  $U_q(\goth g)$ of a
complex Lie algebra $\goth g=\goth s\goth l(n)$ \cite{Drin,Jim,Jimb}.

\subsection{Braid groups and their R-matrix representations}
\lb{subsec2.1}

In the rest of the section we introduce standard notation and recall basic
results on R-matrix representations of the braid groups.

The braid group ${\cal B}_k$ in Artin's presentation
is given by a set of generators $\{\sigma_i\}_{i=1}^{k-1}$
and relations
\ba
\lb{braid1}
\sigma_i\sigma_{i+1}\sigma_i \, =\, \sigma_{i+1}\sigma_i\sigma_{i+1} &&\quad\forall\;
i=1,2,\dots ,k-1,
\\
\lb{braid2}
\sigma_i\sigma_j \, =\, \sigma_j\sigma_i \hspace{10.6mm}&&\quad\forall\; i,j: |i-j|>1\, .
\ea

Let $V$ be a finite dimensional $\mathbb C$-linear space. For any operator
$X\in {\rm End}(V^{\otimes 2})$ and for all integer $i>0$, $j>0$~  denote
\be
\lb{X-k}
X_i\, := I^{\otimes (i-1)}\otimes X\otimes I^{\otimes (j-1)}\in {\rm End}(V^{\otimes (i+j)})\, ,
\ee
where $I\in {\rm Aut}(V)$ is the identity operator.\footnote{Strictly speaking a proper
notation for the l.h.s. of (\ref{X-k}) would be, say, $X_{i}^{(i+j)}$. We use the shortened notation
$X_i$ since a dependence on  $j$
is not critical for our considerations. All formulas below
have sense if the index $j$ is large enough.
A minimal possible value for $j$ in each case
is obvious from the context.}
We also use notation $X_{ij}$ for an operator in ${\rm End}(V^{\otimes k})$,
$1\leq i\neq j\leq k$, acting as $X$ in  {\em component spaces} $V$ with labels $i$ and $j$ and
as identity in the rest. In these notation $X_{i\, i+1}\equiv X_i$.

An operator $\R\in {\rm Aut}(V\otimes V)$
satisfying equality
\be
\lb{ybe}
\R_1 \R_2 \R_1\, =\, \R_2 \R_1 \R_2\, ,
\ee
is called an {\em {\rm R}-matrix}. Any R-matrix generates representations $\Rro$ of the braid groups
${\cal B}_k$, $k=2,3,...$
$$
%\be\lb{R-rep}
\Rro:\quad {\cal B}_k\rightarrow {\rm Aut}(V^{\otimes k}), \quad
\Rro(\sigma_i)= \R_i, \quad 1\leq i\leq k-1.
%\ee
$$
By a slight abuse of notation we assign
the same symbol $\Rro$
to the R-matrix representations of the braid groups ${\cal B}_k$ for different values of index $k$.
This should not cause problems as the braid groups admit a series of monomorphisms
commuting with $\Rro$
\be
\lb{mono-1}
{\cal B}_k\hookrightarrow{\cal B}_{k+1}: \quad\sigma_i\mapsto\sigma_i\quad \forall i=1,\dots k-1.
\ee

\begin{defin}
\lb{defin2.1}
An {\rm R}-matrix  $\R$ is called {\em skew invertible} if there exists an
operator $\Rpsi \in End(V^{\otimes 2})$  such that
\be
\lb{psi}
\Tr{2} \R_{12} \Rpsi_{23} = \Tr{2} \Rpsi_{12} \R_{23} = P_{13} \; .
\ee
Here by $\tr_{\!(i)}$ we denote trace operation in $i$-th space, and by
$P$ --- the permutation operator:
$P(u\otimes v)= v\otimes u\;\; \forall\, u,v\in V$.
\end{defin}

With any skew invertible R-matrix $\R$ we associate a pair of operators
$\D , \C\in {\rm End}(V)$
\be
\lb{dmatr}
{\D}_1 = \Tr{2} \Rpsi_{12} \; , \;\;\; {\C}_2 = \Tr{1} \Rpsi_{12} \; ,
\ee
which, by (\ref{psi}), satisfy equalities
\be
\lb{dmat1}
\Tr{2} \R_{12} {\D}_2  = \id_1 \; , \;\;\; \Tr{1} {\C}_1 \R_{12} = \id_2 \; .
\ee

Further properties of the operators $\D$ and $\C$ are summarized below.

\begin{prop}
\lb{prop2.3} {\rm \cite{I,O}}~
Let $\R$ be a skew-invertible {\rm R}-matrix.
The operators ${\D}$ and ${\C}$ (\ref{dmatr}) satisfy equalities
\ba
\nn
%\lb{dmat2}
\D_1 \, \id_2 = \tr_{\!(3)} \D_3
\R_2^{\pm 1} \, P_{12} \, \R_2^{\mp 1} \, , &&\;
\C_3 \, \id_2 = \tr_{\!(1)} \C_1
\R_1^{\pm 1} \, P_{23} \, \R_1^{\mp 1} \, ,
\\[1mm]
\lb{dmat3}
\R_{12} \, \D_1 \, \D_2 = \D_1 \, \D_2 \, \R_{12}  \, , &&\;
\R_{12} \, \C_1 \, \C_2 = \C_1 \, \C_2 \, \R_{12} \, .
\ea
\end{prop}

Let $W$ be a $\mathbb C$-linear space. For any skew invertible R-matrix $\R$
we define an {\em {\rm R}-trace} map\footnote{This map
is often called a quantum trace or, shortly, a $q$-trace.
In our opinion, the name R-trace is better appropriate to it.
}
$\Rtr: {\rm End}_W(V)\rightarrow W$
\be
\nn
%\lb{R-trace}
Y\mapsto \Rtr(Y) := \tr(\D Y)\, , \quad Y\in {\rm End}_W(V)\, .
\ee

Following properties of the R-trace are simple consequences of the
relations given in the proposition \ref{prop2.3}.

\begin{cor}
\lb{cor2.4}
Let $\R$ be skew invertible {\rm R}-matrix.
For any operator $Y \in {\rm End}_W(V)$
the {\rm R}-trace associated with $\R$
satisfies  relations
\be
\lb{dmat4}
\RTr{2} (\R_{12}^{\varepsilon} \, Y_1 \, \R_{12}^{-\varepsilon}) = I_1\, \Rtr (Y)  \, .
\ee
where $\varepsilon = \pm 1$
and the symbol $\RTr{i}$ denotes taking the {\rm R}-trace in i-th space.

For an element $x^{(k)}\in {\mathbb C}[{\cal B}_k]$
denote $\RX{(k)}:=\rho_R(x^{(k)})\in {\rm End}(V^{\otimes k})$.
Following cyclic property
$$
%\be\lb{dmat5}
\begin{array}{c}
\RTr{1, \dots ,k} \left( \RX{(k)} \ Y^{(k)} \right) =
\RTr{1, \dots ,k} \left(Y^{(k)} \  \RX{(k)} \right)\, .
\end{array}
%\ee
$$
is fulfilled
for any $k\geq 1$ and $Y^{(k)} \in {\rm End}_W(V^{\otimes k})$,
and for all $x^{(k)}\in {\mathbb C}[{\cal B}_k]$.
\end{cor}

\begin{example}
\lb{2.2.3}
{\rm %\small
Permutation $P$: $P (u\otimes v) := v\otimes u\;\;\forall\, u,v\in V$, is the skew invertible R-matrix.
Identity operator $I^{\otimes 2}$ is
the R-matrix which is not skew invertible.
}\end{example}

\begin{example}
\lb{2.2.4}
{\rm %\small
Assume that the quasi-triangular Hopf algebra ${\goth A}_{\cal R}$ admits a representation
$$
%\be\lb{rho-V}
\rho_V: {\goth A}_{\cal R}\rightarrow {\rm End}(V) .
%\ee
$$
As follows from the Yang-Baxter equation (\ref{UYB})
an operator
\be
\lb{R-iz-univ}
\R := \eta P\, (\rho_V \otimes \rho_V)({\cal R}),
\ee
satisfies relation (\ref{ybe}).
Here scaling factor $\eta\in\{\mathbb{C} \setminus 0\}$ is introduced for sake
of future convenience.

The R-matrix (\ref{R-iz-univ}) is skew invertible, its skew inverse matrix is given by formula
(see, e.g., \cite{O}, section 4.1.2)
\be
\nn
\Rpsi\, =\, \eta^{-1} P\,(\rho_V \otimes \rho_V)((id\otimes S){\cal R})
\ee

The matrices $\D$ and $\C$ associated with the R-matrix (\ref{R-iz-univ})
are:
\be
\lb{CD-iz-univ}
\D\, =\, \eta^{-1}\,\rho_V(u)\, , \qquad
\C\, =\, \eta^{-1}\,\rho_V(S(u))\, .
\ee
Both, they are invertible and their properties (\ref{dmat3}) are descending from
(\ref{Ruu}).
}\end{example}

\subsection{Hecke algebras and Hecke type R-matrix}
\lb{subsec2.3}

An A-type {\em Hecke algebra} ${\cal H}_k(q)$ is a quotient algebra of the group algebra
${\mathbb C}[{\cal B}_k]$ \eqref{braid1}, \eqref{braid2} by relations
$$
%\be\lb{Hecke}
(\sigma_i -q 1)(\sigma_i +q^{-1} 1)=0\qquad \forall\, 1\leq i\leq k-1.
%\ee
$$
Under following conditions on the  parameter $q$
\be
\lb{restrict}
\hspace{-29mm}\mbox{\rm\bf [k]} \hspace{37mm}
i_q := {(q^i-q^{-i})/(q-q^{-1})}\neq 0 \quad \forall i=2,3,\dots ,k,\hspace{5mm}
\ee
the algebra ${\cal H}_k(q)$  is
isomorphic to the group algebra of the symmetric group ${\mathbb C}[S_k]$ and, hence, semisimple.
It's irreducible representations as well as its central idempotents are labeled by
a set of partitions $\lambda\vdash k$.
We are particularly interested in a series of idempotents
corresponding to the one dimensional representations $\lambda = (1^k),\, k=1,2,\dots$ .
These idempotents -- we denote them as $a^{(k)}$ -- admit a recursive construction
(see, e.g., \cite{HIOPT}, section~1, or \cite{GPS}, section~2.3, or \cite{TW}, lemma 7.2)
\ba
\lb{q-anti1}
\hspace{-5mm}
a^{(1)} = 1, &&
a^{(k)} = \frac{(k-1)_q}{k_q}\,
a^{(k-1)}\,
\Bigl(\frac{q^{k-1}}{(k-1)_q}\, 1\, -\, \sigma_{k-1}\Bigr)
a^{(k-1)}
\\[1mm]
\lb{q-anti2}
&& =\,\frac{(k-1)_q}{k_q}\,
a^{(k-1)\uparrow 1}\,
\Bigl(\frac{q^{k-1}}{(k-1)_q}\, 1\, -\, \sigma_1\Bigr)
a^{(k-1)\uparrow 1}\qquad \forall\, k=2,3,\dots \, ,
\ea
where  we use symbol $x^{(k)\uparrow 1}\in {\cal H}_{k+1}(q)$
for an image of  element $x^{(k)}\in {\cal H}_k(q)$ under following
algebra  monomorphism (c.f. with (\ref{mono-1})):
$$
{\cal H}_k\hookrightarrow{\cal H}_{k+1}: \quad\sigma_i\mapsto
\sigma_{i+1}\quad \forall i=1,\dots k-1\, .
$$
The idempotents $a^{(k)}$
obey relations
\ba
\lb{idemp-1}
a^{(k)}\sigma_i  =\sigma_i a^{(k)} = -q^{-1} a^{(k)}&&
\qquad\;\;\forall\; i=1,2,\dots ,k-1\, ,
\\[2mm]
\lb{idemp-2}
a^{(k)} a^{(i)\uparrow j} = a^{(i)\uparrow j} a^{(k)} = a^{(k)}\, , &&
\qquad\;\;\, \mbox{if~~}\; i+j\leq k\, .
\ea

\vspace{1mm}
An R-matrix $\R$ satisfying quadratic minimal characteristic identity
is called a {\em Hecke type} R-matrix. By  an appropriate  rescaling of $\R$
one always can turn its characteristic identity
to a form
\be
\lb{H}
(\R -q I)(\R+q^{-1}I) = 0\, .
\ee
In this case the corresponding representations $\rho_R$
become representations of the Hecke algebras ${\cal H}_k(q)$
\be
\lb{R-rep2}
\Rro:\quad {\cal H}_k(q)\rightarrow {\rm Aut}(V^{\otimes k}), \quad
\Rro(\sigma_i)= \R_i, \quad 1\leq i\leq k-1.
\ee
We reserve special notation for the R-matrix images of idempotents $a^{(k)}$:
\be
\lb{A-k}
A^{(k)}:=\Rro(a^{(k)})\, , \quad
A^{(k)\uparrow 1}:=\Rro(a^{(k)\uparrow 1})\quad \forall\, k\geq 1\, .
\ee
We also put $A^{(0)}:=1$. The elements
$A^{(k)}$ will be further referred as {\em $k$-antisymmetrizers}.

\begin{rem}{\rm\small
The R-matrix analogues of relations (\ref{q-anti1})--(\ref{idemp-2})
have been described in literature (see \cite{Jimb}, \cite{Gur})
even earlier then their algebraic prototypes.
}
\end{rem}

\subsection{$GL_q(n)$ type R-matrix}
\lb{subsec2.4}

\begin{defin}
\lb{def2.4}
Consider a Hecke type {\rm R}-matrix $\R$. Assume that parameter $q$
in its characteristic identity (\ref{H}) satisfies conditions
{\rm\bf  [n]} (\ref{restrict}), so that
antisymmetrizers $A^{(2)}, \dots ,A^{(n)}$ are well defined.
$\R$ is called {\em $GL_q(n)$ type {\rm R}-matrix} if two conditions
\ba
\lb{GL(n)}
A^{(n)} \Bigl( {q^n\over n_q}\,I\,-\,\R_n\Bigr)A^{(n)}&=& 0\hspace{6cm}
\\[-2mm]
\nn
\mbox{and}\hspace{77mm}&&
\\[-2mm]
\lb{GL(n)-2}
\qquad {\rm rk}A^{(n)} &=& 1
\ea
are fulfilled.
\end{defin}
\begin{rem}
{\rm\small
Assuming $(n+1)_q\neq 0$,  the condition (\ref{GL(n)})
is equivalent to $A^{(n+1)}=0$.
For generic values of $q$, assuming
validity of (\ref{GL(n)}),
the condition (\ref{GL(n)-2})
is equivalent to demanding skew invertibility of $\R$
(see \cite{Gur}, propositions 3.6 and 3.10).
}
\end{rem}

\begin{prop}\lb{prop2.9} {\rm \cite{Gur,I}}~
Let $\R$ be a skew invertible {\rm R}-matrix of the type $GL_q(n)$. Then
$\C$ and $\D$ are invertible and
following relations are fulfilled
\ba
\lb{C-D}
\D\, \C& =& \C\, \D\, =\, q^{-2n}\, I ,
\\[1mm]
\lb{A-A}
\RTr{k} A^{(k)}& =& q^{-n}\, {(n+1-k)_q\over k_q} A^{(k-1)}
\quad \forall\, k=1,2,\dots ,n,
\\[-1mm]
\lb{D-A}
A^{(n)}\prod_{i=1}^n(\D)_i &=&
\prod_{i=1}^n(\D)_i \, A^{(n)}\, =\, q^{-n^2}A^{(n)} .
\ea
\end{prop}

\begin{example}
\lb{2.4.4}
{\rm %\small
Consider the case ${\goth A}_{\cal R}$ is
the quantized universal enveloping algebra $U_q\goth s\goth l(n)$.
Let $V$ be a vector representation of $U_q\goth s\goth l(n)$,
$\dim V =n$.
In this case formula (\ref{R-iz-univ}) with the scaling factor
chosen as $\eta = q^{1/n}$
gives a standard {\em Drinfeld-Jimbo's}
 R-matrix $\R^\circ$ of the $GL_q(n)$ type (see \cite{KSch}, section 8.4.2):
\be
\lb{DJ}
\R^{\circ}\, =\,
 \sum_{i,j=1}^n q^{\delta_{ij}}E_{ij}\otimes E_{ji}\, +\,
(q-q^{-1})\,\sum_{i<j } E_{ii}\otimes E_{jj}\, .
\ee
Here $(E_{ij})_{kl}:=\delta_{ik}\delta_{jl}$, $\, i,j=1,\dots ,n,\,$
is a standard basis of $n\times n$ matrix units.
Via the so-called twist procedure (for details see \cite{Resh2})
$\R^{\circ}$ gives rise to a multiparametric family of $GL_q(n)$ type R-matrices
\be
\lb{DJ-multi}
\R^{\,f}\, :=\, F \R^{\circ} F^{-1}\, =\,
 \sum_{i,j=1}^n q^{\delta_{ij}}{f_{ij}\over f_{ji}}\, E_{ij}\otimes E_{ji}\, +\,
(q-q^{-1})\,\sum_{i<j } E_{ii}\otimes E_{jj}\, ,\quad \forall\, f_{ij}\in \{\mathbb{C}\setminus 0\}.
\ee
Here $F:=\sum_{i,j=1}^n f_{ij}\, E_{ii}\otimes E_{jj}\,$ is a {\em twisting R-matrix}.
In what follows we use these particular R-matrices for illustration purposes.
Their corresponding matrices $D_{\R^\circ}$ and $D_{\R^f}$ are
\be
\nn
D_{\R^\circ}\, =\, D_{\R^f}\, =\, \sum_{i=1}^n q^{2(i-n)-1} E_{ii}\, .
\ee
}\end{example}

\begin{rem}
\lb{2.4.5}
{\rm\small
Generally speaking,  $GL_q(n)$ type R-matrix can be realized
in a tensor square of space $V$ whose dimension is different from $n$.
Examples of the R-matrices
for any $\dim V\geq n$
are given in \cite{Gur}, in section 4.
In what follows we do not assume any relation between the parameter
$n$ in the definition \ref{def2.4} and the dimension of the space $V$, unless it is stated explicitly.
}\end{rem}

\section{Quantized functions on a cotangent bundle over matrix  group}
\lb{sec3}

In this section we recall
definition of a quantum group cotangent bundle and develop
in linear cases -- $GL_q(n)$ and $SL_q(n)$ --
a basic techniques for it's structure investigation.

\subsection{Quantized functions over matrix group (RTT algebra)}
\lb{subsec3.1}

\begin{defin} {\rm \cite{Drin,FRT}}
\lb{def3.1}
Let $\R$
be a skew invertible  {\rm R}-matrix.
An associative unital algebra   generated by a set of matrix components
${\cal k}T^i_j{\cal k}_{i,j=1}^{\dim V}$ satisfying relations
\be
\lb{RTT}
\R_{12} \, T_1 \, T_2 = T_1 \, T_2 \, \R_{12}
\ee
is  denoted as ${\goth F}[\R]$ and called an {\em RTT algebra}.
The RTT algebra is endowed in a standard way with
the coproduct and the counit
\be
\lb{bi-RTT}
\Delta (T^i_j) =  \sum_k T^i_k\otimes T^k_j\, , \qquad \epsilon(T^i_j)= \delta^i_j\, .
\ee
Let further extend the RTT algebra by a set of inverse matrix components
${\cal k}(T^{-1})^i_j{\cal k}_{i,j=1}^{\dim V}$:
\be
\lb{T-inverse}
\sum_k T^i_k\, (T^{-1})^k_j\, =\, \sum_k (T^{-1})^i_k\, T^k_j\, =\,
\delta^i_j\, 1\, .
\ee
The extended algebra can be endowed with the antipode mapping
\be
\lb{antipod}
S(T^i_j)\, =\, (T^{-1})^i_j\, , \qquad\mbox{so that (see \cite{Resh}):}\qquad
S^2(T)\,\D\, =\, \D\, T\, .
\ee
The resulting Hopf algebra is further denoted as ${\goth F\goth G}[\R]$.
\end{defin}

\begin{example}
\lb{3.1.2}
{\rm %\small
Consider the quasi-triangular Hopf algebra ${\goth A}_{\cal R}$ together with its representation
$\rho_V$ (see example \ref{2.2.4}).
For any $x\in{\goth A}_{\cal R}$
denote ${\cal k}\rho_V(x)^i_j{\cal k}$
a matrix of the operator $\rho_V(x)$
in a certain basis in the space $V$.

Let ${\goth A}_{\cal R}^*$ be the dual Hopf algebra and let
${\cal h}\cdot ,\cdot {\cal i}$ denote a non degenerate pairing between ${\goth A}_{\cal R}$ and
${\goth A}_{\cal R}^*$. Consider two matrices of linear functionals  on ${\goth A}_{\cal R}$ ---
$\,T^i_j$ and $(T^{-1})^i_j\,$ ---\, such that
\be
\lb{motiv1}
\langle T^i_j,x\rangle = \rho_V(x)^i_j\, ,\quad
\langle(T^{-1})^i_j\, , x\rangle  = \rho_V(S(x))^i_j\quad \forall\, x\in {\goth A}_{\cal R}\, .
\ee
It is easy to see that these functionals satisfy conditions of the
definition \ref{def3.1} (for details see, e.g., \cite{Bur}),
the numeric R-matrix $\R$ in (\ref{RTT}) in this case is given by (\ref{R-iz-univ}),
relation (\ref{antipod}) for the square of antipode descends from (\ref{square-antipode}).
The functionals $T^i_j$ and $(T^{-1})^i_j$ generate a Hopf subalgebra in
${\goth A}_{\cal R}^*$.

In case if ${\goth A}_{\cal R}$ is a universal enveloping algebra $U{\goth g}$ of some Lie algebra $\goth g$,
the dual Hopf algebra $(U\goth g)^*$ can be treated as
${\rm Fun}(\goth G)\equiv {\goth F\goth G}$,
where $\goth G$ is a formal group corresponding to $\goth g$. Therefore, heuristically
we can treat the RTT algebras ${\goth F\goth G}[\R]$ and ${\goth F}[\R]$
as algebras of {\em quantized functions over matrix group and matrix semigroup},
respectively.
Here the term {\em matrix} refers to a matrix form of the coproduct (\ref{bi-RTT}); the term
{\em quantized} means that relations (\ref{RTT}) in general define a noncommutative
product.
}\end{example}

In the rest of the subsection we describe a construction of the inverse matrix
$T^{-1}$ for the RTT algebra associated with the $GL_q(n)$ type R-matrix.
\medskip

Consider an element
\be
\lb{detT2}
\Rdet T := \Tr{1, \dots ,n} \bigl( A^{(n)} T_1 T_2 \dots T_n \bigr) .
\ee
By the definition of the coproduct \eqref{bi-RTT}  and due to the rank 1 condition (\ref{GL(n)-2})
the element $\Rdet T$ is group-like
\be
\nn
%\lb{grouplike}
\Delta (\Rdet T) = \Rdet T\otimes\Rdet T\, ,
\ee
and it satisfies relations
\be
\nn
%\lb{detT}
A^{(n)} \, T_1 T_2 \dots T_n = T_1 T_2 \dots T_n \,
A^{(n)} = A^{(n)} \, \Rdet T \, .
\ee
Therefore, it is natural to call ${\det}_R T$  a {\em  determinant
of the matrix  $T$}.

\begin{prop} {\rm \cite{Gur}}
\lb{3.1.3}
Let $\R$ be a skew invertible $GL_q(n)$ type {\rm R}-matrix.
The following relation is satisfied in
the corresponding RTT algebra ${\goth F}[\R]$
\be
\nn
%\lb{det}
(\Rdet T)\, T\,  =\,  (\O T\O^{\!\!\!-1})\, \Rdet T\, ,
\ee
where $~\O,\, \O^{\!\!\!-1}\in {\rm Aut}(V)$ are mutually inverse matrices:
\ba
\lb{O}
\O_1 &=& n_q\, \Tr{2,\dots ,n+\! 1}\left(P_1 P_2\dots P_n\, A^{(n)}\right)\, ,
\\[2mm]
\nn
%\lb{O-inv}
(\O^{\!\!\!-1})_1 &=& n_q\, \Tr{2,\dots ,n+\! 1}\left(A^{(n)}P_n\dots  P_2 P_1 \right)\, ,
\ea
(recall that $P_i$ are permutation operators acting in components spaces $V_i\otimes V_{i+1}$).
\end{prop}

\begin{cor}
\lb{3.1.3-cor}
In the assumptions of proposition \ref{3.1.3} consider an extension
of the RTT algebra ${\goth F}[\R]$ by an element $(\Rdet T)^{-1}$ subject to relations
\be
\nn
%\lb{det-inv}
(\Rdet T)^{-1}\, T = (\O^{\!\!\!-1}T\O)(\Rdet T)^{-1} \, ,\quad
\Rdet T\, (\Rdet T)^{-1} = (\Rdet T)^{-1}\, \Rdet T = 1\, ,
\ee
In the extended algebra the inverse matrix $T^{-1}$ satisfying relations
(\ref{T-inverse}) is given by formula
\be
\nn
%\lb{T-inv}
(T^{-1})_1\, =\,
q^{n(n-1)} n_q\, \RTr{2,\dots ,n}\left(T_2\dots T_n A^{(n)}\right) (\Rdet T)^{-1}\, ,
\ee
The resulting Hopf algebra is called {\em $GL_q(n)$ type RTT algebra} and denoted as
${\goth F {\scriptstyle GL_q(n)}}[\R]$.

Assume additionally that for the {\rm R}-matrix $R$ the corresponding
matrix $\O$ (\ref{O}) is scalar: $\O\propto I$.
In this case $R$ is called  the R-matrix of {\em $SL_q(n)$ type}.
In the corresponding RTT algebra
${\goth F {\scriptstyle GL_q(n)}}_R$
the element $\Rdet T$ is central. A quotient of this algebra by relation
$\Rdet T\, =\, 1$ is called {\em $SL_q(n)$ type RTT algebra}
and denoted as ${\goth F {\scriptstyle SL_q(n)}}[\R]$.
\end{cor}

\begin{rem}
{\rm\small  For a skew invertible $GL_q(n)$ type R-matrix $\R$
consider a system of equations
$$
\R_{12}\, N_1\, N_2\, =\, N_1\, N_2\, \R_{12}\, , \qquad
N^n\, \propto\, \O\qquad \mbox{for some $N\in {\rm Aut}(V)$.}
$$
Note that a consistency condition for these equations\, ---
$\,\R_{12}\O_1 \O_2 = \O_1 \O_2 \R_{12}\,$
---\, is satisfied (see \cite{OP}). By any solution $N$ of these equations one can construct
the $SL_q(n)$ type R-matrix
\be
\lb{partner-R}
\widetilde\R_{12}\, :=\, N_1 \R_{12} N_1^{-1}\, =\, N_2^{-1}\R_{12} N_2\, .
\ee
}\end{rem}
\begin{example}{\rm %\small
For the R-matrices  described in the example \ref{2.4.4}
one has
\be
\nn
O_{\R^\circ}\, =\, - I\, , \qquad
O_{\R^f}\, =\, -\,
{\textstyle \sum_{i=1}^n}\Bigl({\textstyle \prod_{j\neq i}}{f_{ji}/ f_{ij}}\Bigr)\, E_{ii}\, .
\ee
So, $\R^\circ$ is $SL_q(n)$ type, while $\R^f$ is $SL_q(n)$ type only if
$\forall\, i=1,\dots ,n:\;\prod_{j\neq i}(f_{ji}/f_{ij})=
{}^{\raisebox{1.6pt}{$\scriptstyle n$}}\hspace{-4.2pt}\sqrt{1}$.
Taking a diagonal $n$-th root $O_{\R^f}^{1/n}$ of the diagonal matrix $O_{\R^f}$ one finds the
$SL_q(n)$ type R-matrix associated with $\R^f$:
\be
\nn
\widetilde\R^f\, =\, \R^{\tilde f}\, , \quad \mbox{where}\quad
{\tilde f}_{ij}\, :=\,{\textstyle \prod_{k\neq i,j}} (f_{ij} f_{jk} f_{ki})^{1/n}\, , \quad
\mbox{so that}\quad O_{\R^{\tilde f}} = - I\, .
\ee
}\end{example}

\subsection{Quantized right invariant vector fields (reflection equation algebra)}
\lb{subsec3.2}

\begin{defin} {\rm \cite{KS}}
\lb{3.2.1}
Let $\R$
be a skew invertible  {\rm R}-matrix.
An  associative unital algebra ${\goth L\goth G}[\R]$
generated by a set of matrix components
${\cal k}L^i_j{\cal k}_{i,j=1}^{\dim V}$ satisfying relations
\be
\lb{REA}
L_1 \R_{12} L_1 \R_{12}\, =\, \R_{12} L_1 \R_{12} L_1
\ee
is called a {\em reflection equation algebra} or, shortly, {\em RE algebra}.
The RE algebra ${\goth L\goth G}[\R]$ is naturally endowed with a structure
of left coadjoint ${\goth F\goth G}[\R]$-comodule algebra
\be
\lb{REA-com}
\delta_{\ell}(L^i_j) = \sum_{k,m}\, T^i_k (T^{-1})^m_j \otimes L^k_m\, .
\ee
\end{defin}

\begin{example} {\rm \cite{FRT}}
\lb{3.2.2}
{\rm %\small
In notations of the examples  \ref{2.2.4}, \ref{3.1.2}
consider following ${\goth A}_{\cal R}$-valued matrices
\be
\lb{dop1}
\begin{array}{rcl}
{L^{\scriptscriptstyle (+)}}^i_j\, =\, \langle id \otimes T^i_j, \,
{\cal R} \rangle\, ,&&
{L^{\scriptscriptstyle (-)}}^i_j\, =\,  \langle S(T^i_j) \otimes id\, ,
{\cal R} \rangle\, =\,
\langle T^i_j \otimes id\, ,  {\cal R}^{-1} \rangle\, ,
\\[2mm]
((L^{\scriptscriptstyle (+)})^{-1})^i_j =
\langle id \otimes T^i_j, \, {\cal R}^{-1} \rangle ,&&
((L^{\scriptscriptstyle (-)})^{-1})^i_j =
\langle   T^i_j\otimes id\, , \, {\cal R}\rangle .
\end{array}
\ee
As a consequence of the Yang-Baxter equation (\ref{UYB})  components
of these matrices satisfy relations
\be
\lb{dop2a}
\begin{array}{c}
\R_{12} \, L^{\scriptscriptstyle( \pm)}_2 \, L^{\scriptscriptstyle (\pm)}_1 =
L^{\scriptscriptstyle (\pm)}_2 \, L^{\scriptscriptstyle (\pm)}_1 \, \R_{12}\,  , \quad
\R_{12} \, L^{\scriptscriptstyle (+)}_2 \, L^{\scriptscriptstyle (-)}_1 =
L^{\scriptscriptstyle (-)}_2 \, L^{\scriptscriptstyle (+)}_1 \, \R_{12} \, ,
\end{array}
\ee
where $\R$ is given by (\ref{R-iz-univ}).
By (\ref{coco}), the elements $\bigl({(L^{\scriptscriptstyle (\pm)})^{\pm 1}}\bigr)^i_j$ generate a
Hopf ${\goth A}_{\cal R}$-subalgebra
$$
%\be\lb{bialg-L}
\Delta({L^{\scriptscriptstyle (\pm)}}^i_j) =
\sum_k{L^{\scriptscriptstyle (\pm)}}^i_k \otimes {L^{\scriptscriptstyle (\pm)}}^k_j\, , \quad
\epsilon({L^{\scriptscriptstyle (\pm)}}^i_j) = \delta^i_j \, ,\quad
S({L^{\scriptscriptstyle (\pm)}}^i_j) =
((L^{\scriptscriptstyle (\pm)})^{-1})^i_j\, .
%\ee
$$
Consider a composite matrix $L$ with components
\be
\lb{L-iz-univ}
L^i_j\, :=\, q^{(n-{1\over n})}\sum_k\, ((L^{\scriptscriptstyle (-)})^{-1})^i_k
{L^{\scriptscriptstyle (+)}}^k_j\, =\,q^{(n-{1\over n})}
\langle  id\otimes T^i_j\, ,
{\cal R}_{21}{\cal R}_{12}\rangle\, ,
\ee
where our choice of a numeric factor $q^{n-{1\over n}}$
is argued in appendix \ref{append-A}.
By (\ref{dop2a}), components of  $L$ satisfy reflection equation
 (\ref{REA}), where $\R$ is given by
(\ref{R-iz-univ}). Note that  an ${\goth A}_{\cal R}$-subalgebra generated by
the elements $L^i_j$ (\ref{L-iz-univ}) does not carry a natural Hopf algebra structure.
Instead, it obeys a coadjoint comodule algebra structure (\ref{REA-com}) with respect to the Hopf
${\goth A}^*_{\cal R}$-subalgebra generated by the components of the matrices $T$
and $T^{-1}$ (\ref{motiv1}).
}
\end{example}

Let us comment on a geometric interpretation of the RE algebra.
In \cite{FRT}
the matrices $L^{(\pm)}$ were used to develop an RTT type description for
the quantized universal enveloping algebra $U_q \goth g$.
Consider the case $\goth g={\goth s\goth l}(n)$ and let $V$ be its vector representation.
The corresponding  $GL_q(n)$ type R-matrix $\R$ is given
in the example \ref{2.4.4}.
Making a linear change of generators $L^i_j \rightarrow \ell^i_j$:
\be
\lb{ell}
L^i_j\, =\, \delta^i_j + (q-q^{-1}) \ell^{\,i}_j\, .
\ee
and using the Hecke condition (\ref{H})
the reflection equation (\ref{REA}), for $q^2\neq 1$,  can be equivalently rewritten as
\be
\lb{REA2}
\ell_1 \R_{12}\ell_1 \R_{12} \, -\, \R_{12}\ell_1 \R_{12}\ell_1\,
=\, \R_{12}\ell_1 \, -\, \ell_1 \R_{12}\, .
\ee
In a "classical" limit $q\rightarrow 1$ the R-matrix (\ref{DJ}) tends to the permutation and the
equations (\ref{REA2}) go into commutation relations for the basis
of generators of the
Lie algebra
$\goth g\goth l(n)$
\be
\lb{REA-limit}
[\ell_1, \ell_2]\, =\, P_{12}(\ell_1 - \ell_2)\, .
\ee
Classically we can treat ${\cal k}\ell^i_j{\cal k}_{i,j=1}^{p}$ as
a basis  of  right invariant vector fields on $GL(n)$.
Transformation of these basic fields under the left
transition by a group element $t\in GL(n)$ is given by formula (c.f. with (\ref{REA-com}))
$$
\delta_{\ell}(t)\, :\;\; \ell^i_j\mapsto  \sum_{k,m=1}^n\, t^i_k\, \ell^k_m\, (t^{-1})^m_j\, ,
\qquad\mbox{where}\quad
t^i_j:=
%\langle T^i_j,t\rangle =
\rho_V(t)^i_j\, .
$$
Extrapolating this interpretation to a "quantum" case $q\neq 1$ we
call ${\cal k}L^i_j{\cal k}_{i,j=1}^{n}$
{\em a  basis of quantized right invariant vector
fields over matrix group}.
\medskip

It is technically convenient to
introduce notation
\ba
\lb{L-over}
L_{\overline{1}}:= L_1\, , &&
L_{\overline{k+1}} := \R_k \, L_{\overline{k}} \, \R_k^{-1}\, ,
\\[2pt]
\nn
%\lb{L-under}
L_{\underline{1}}:= L_1\, , &&
L_{\underline{k+1}} := \R_k^{-1} \, L_{\underline{k}} \, \R_k\qquad \forall\, k\geq 1 .
\ea
In terms of these {\em {\rm R}-copies} $L_{\overline k}$, $L_{\underline k}$ of the matrix $L$
the reflection equation (\ref{REA}) can be equivalently written in any of the following forms
\be
\lb{REA-forms}
R_k\, L_{\overline{k}}\,L_{\overline{k+1}}\, =\, L_{\overline{k}}\,L_{\overline{k+1}}\, R_k\, ,
\qquad
R_k\, L_{\underline{k+1}}\,L_{\underline{k}}\, =\, L_{\underline{k+1}}\,L_{\underline{k}}\, R_k\qquad
\forall\, k\geq 1\, .
\ee
Taking into account commutativity relations
\be
\lb{LR-comm}
R_i\, L_{\overline k}\, =\,L_{\overline k}\, R_i\, , \qquad
R_i\, L_{\underline k}\, =\,L_{\underline k}\, R_i\qquad \forall\, i,k:\; k\neq i,i+1\, ,
\ee
one sees that the R-copies $L_{\overline k}~(L_{\underline k})$
 of the matrix $L$ in the RE algebra ${\goth L\goth G}[\R]$ formally satisfy the same
relations as the usual copies $T_k~( T_k^{-1})$ of the matrix $T~(T^{-1})$ in the RTT algebra
${\goth F\goth G}[\R]$.

Matrix monomials in two different series of the R-copies
satisfy relations
\be
\lb{lll}
L_{\overline{1}} \, L_{\overline{2}} \dots L_{\overline{k}} =
L_{\underline{k}} \dots L_{\underline{2}} \, L_{\underline{1}}\qquad \forall\, k\geq 1\, ,
\ee
For $k=2$ the  equality (\ref{lll}) is identical to the
reflection equation (\ref{REA}). For $k>2$
this equality follows by induction on $k$.
Note that monomials (\ref{lll}) transform covariantly
under the left coadjoint coaction (\ref{REA-com})
\be
\lb{lll-trans}
\delta_\ell\left(L_{\overline{1}}  \dots L_{\overline{k}}\right)\, =\,
\left(T_1\dots T_k\otimes 1)(1\otimes L_{\overline{1}}  \dots L_{\overline{k}}\right)
\left(S(T_1\dots T_k)\otimes 1\right) .
\ee

\medskip
Following proposition goes back to theorem 14  from \cite{FRT}
(see also \cite{I}, proposition 5).

\begin{prop} Let $\R$
be a skew invertible  {\rm R}-matrix.
For an element $x^{(k)}\in {\mathbb C}[{\cal B}_k]$
denote
\be
\lb{center}
{ch}(x^{(k)}) := \RTr{1\dots k} \left(\RX{(k)} \,
L_{\overline{1}} \, L_{\overline{2}} \dots L_{\overline{k}}  \right)
% \, =\, \RTr{1\dots k} \left(
% L_{\overline{1}} \, L_{\overline{2}} \dots L_{\overline{k}}\, \RX{(k)}  \right)
,
\ee
where $\RX{(k)}:=\rho_R(x^{(k)})\in {\rm End}(V^{\otimes k})$.
Consider a linear subspace ${\goth C\goth h}[\R]\subset {\goth L\goth G}[\R]$ spanned by the
unity and by elements $ch(x^{(k)})\;\; \forall k\geq 1$ and
$\forall x^{(k)}\in {\mathbb C}[{\cal B}_k]$.
The space
${\goth C\goth h}[\R]$ is a subalgebra of the
center of the RE algebra ${\goth L\goth G}[\R]$.
It is called {\em a characteristic subalgebra} of the RE algebra ${\goth L\goth G}[\R]$.
The characteristic subalgebra is invariant with respect
to the left ${\goth F\goth G}[\R]$ coadjoint coaction (\ref{REA-com}).
\end{prop}

\ni
{\em Proof.}
In a setting of the quasi-triangular Hopf algebras these statements were proved in
\cite{Drinf,Resh} (see there sec.3 and sec.4, respectively).
Below we prove the proposition in the RE algebra setting.

\medskip
Consider an arbitrary element $ch(x^{(k)})$ of the characteristic subalgebra.
We first prove a following version of the formula (\ref{center})
\ba
\nn
ch(x^{(k)})\,  I_1 &=& \RTr{2, \dots ,k+1} \left(\RX{(k)\uparrow 1} \,
L_{\overline{2}} \, L_{\overline{3}} \dots L_{\overline{k+1}}  \right) =\hspace{5mm}
\\[2mm]
\lb{cent1}
&=& \RTr{2, \dots ,k+1} \left(\RX{(k)\uparrow 1} \,
L_{\underline{k+1}} \dots  L_{\underline{3}} \, L_{\underline{2}}  \right)  .
\ea
Here the first equality results from a calculation
\ba
\nn
\lefteqn{\RTr{2, \dots ,k+1} \left(\RX{(k)\uparrow 1} \,
L_{\overline{2}} \dots L_{\overline{k+1}}  \right)}\hspace{15mm}&&
\\[2mm]
\nn
&=&\RTr{2, \dots ,k+1} \left(\RX{(k)\uparrow 1} \, \R_1 \cdots \R_k \,
L_{\overline{1}} \dots L_{\overline{k}} \, \R_k^{-1} \cdots \R_1^{-1} \right)
\\[2mm]
\nn
&=& \RTr{2, \dots ,k+1} \left(  \R_1 \cdots \R_k \, (\RX{(k)} \,
L_{\overline{1}} \dots L_{\overline{k}}) \, \R_k^{-1} \cdots \R_1^{-1} \right)
\\[2mm]
% \nn
% &=& \RTr{2, \dots ,k} \left\{  \R_1 \cdots \R_{k-1} \,  \left(I_k \RTr{k'} (\RX{(k)} \,
% L_{\overline{1}} \dots L_{\overline{k'}})\right) \, \R_{mk-1}^{-1} \cdots \R_1^{-1} \right\}
% \\[2mm]
\nn
&=& \dots\;\, =\;\,  \RTr{1, \dots ,k} \left(\RX{(k)} \,
L_{\overline{1}} \, L_{\overline{2}} \dots L_{\overline{k}}  \right)  ,
\ea
where in the last line we applied (\ref{dmat4}) $k$ times.
To prove
the second equality in (\ref{cent1}) we first use the relation (\ref{lll})
and then perform similar transformations.

With the use of (\ref{cent1})  and (\ref{lll}) checking centrality of  $ch(x^{(k)})$
is straightforward
\ba
\nn
L_1 \, ch(x^{(k)}) &=& \RTr{2, \dots ,k+1} \left(\RX{(k)\uparrow 1} \,
L_{\overline{1}} L_{\overline{2}} \, L_{\overline{3}} \dots L_{\overline{k+1}}  \right)
\\[2mm]
\nn
&=& \RTr{2, \dots ,k+1} \left(\RX{(k)\uparrow 1} \,
L_{\overline{k+1}} \dots L_{\overline{2}} \,
L_{\overline{1}} \right) = ch(x^{(k)}) \, L_1 \, .
\ea

The invariance of $ch(x^{(k)})$ under
the left ${\goth F\goth G}[\R]$ coadjoint coaction follows immediately from
(\ref{lll-trans}) together with the relation (\ref{antipod}) for the square of antipode.
\hfill$\blacksquare$

\medskip
Consider a series of elements of the RE algebra ${\goth L\goth G}[\R]$
\be
\lb{s}
p_i \, := \,  \Rtr (L^i) \, ,\quad i=1,2,\dots\, .
\ee
Further on they are called {\em power sums}.
Following calculation
$$
L_1\, p_i\, =\, \RTr{2} L_1 \R_{12} L_1^i \R_{12}^{-1}\, =\, \RTr{2}
\R_{12}^{-1} L_1^i \R_{12} L_1\, =\,p_i\, L_1\, ,
$$
proves centrality of the power sums. Here in the first and the last equalities we use
formula (\ref{dmat4}), and
the second equality is a consequence of (\ref{REA}).
Actually, the power sums belong to the characteristic subalgebra
${\goth C\goth h}[\R]$:
$$
p_i\, =\, ch(\sigma_{i-1}\dots \sigma_2\,\sigma_1)\, ,
$$
which is verified by a following transformation
\ba
\nn
ch(\sigma_{i-1}\dots \sigma_2\,\sigma_1)
&=& \RTr{1, \dots ,i} \left(L_{\overline{1}} \dots L_{\overline{i}}\,
(\R_{i-1} \dots  \R_1) \right)
\\[2mm]
\nn
&&
\hspace{-15mm}
=\;\, \RTr{1, \dots ,i} \left(
L_{\overline{1}} \dots L_{\overline{i-1}} \, (\R_{i-1} \dots \R_1) L_{1}
(\R_1^{-1} \dots \R_{i-1}^{-1}) \, (\R_{i-1} \dots  \R_1) \right)
\\[2mm]
\nn
&&
\hspace{-15mm}
=\;\, \RTr{1, \dots ,i-1} \left( L_{\overline{1}} \dots L_{\overline{i-1}} \,
\left(\RTr{i} \R_{i-1}\right) (\R_{i-2} \dots \R_1) L_{1} \right)
\\[2mm]
\nn
&&
\hspace{-15mm}
=\;\, \RTr{1, \dots ,i-2} \left( L_{\overline{1}} \dots L_{\overline{i-2}} \,
\left(\RTr{i-1} \R_{i-2}\right) (\R_{i-3} \dots \R_1) L_{1}^2 \right)
\\[2mm]
\nn
&&
\hspace{-15mm}
=\;\, \dots\;\, =\;\,  \Rtr (L^i)\;\, =\;\, p_i \, .
\ea
Here  we repeatedly expand the notation
$L_{\overline{j}}=(\R_{j-1} \dots \R_1) L_{1}
(\R_1^{-1} \dots \R_{j-1}^{-1})$ for $j=i,\dots ,2$,
and use (\ref{dmat1}).

\medskip
Let $\R$ be a skew invertible R-matrix of the Hecke type.
Assuming that conditions {\rm\bf  [k]} (\ref{restrict})
are fulfilled consider a series of elements
$a_i\in {\goth C\goth h}[\R]$, $i=0,1,\dots\, k,$
in the corresponding {\em Hecke type} RE algebra ${\goth L\goth G}[\R]$
\be
\lb{sigma}
a_0\, :=\, 1\, , \qquad
a_i \, := ch(a^{(i)})\, =
\, \RTr{1, \dots ,i} \left( A^{(i)}
L_{\overline{1}} \dots L_{\overline{i}} \right)\qquad  \forall\, 1\leq i\leq k\, ,
\ee
where notations $a^{(i)}$, $A^{(i)}$ were explained in (\ref{q-anti1}), (\ref{A-k}).
The elements $a_i$ are called
{\em elementary symmetric functions}.

\begin{defin}
\lb{3.2.4}
Let $\R$ be a skew invertible $GL_q(n)$ type {\rm R}-matrix.
A central extension of the corresponding RE algebra ${\goth L\goth G}[\R]$
by an element $a_n^{-1}:\; a_n a_n^{-1}\, = \, 1$ is called {\em $GL_q(n)$
type RE algebra} and denoted as
${\goth L\scriptstyle GL_q(n)}[\R]$.
A quotient of this algebra by a relation
\be
\lb{RE-sl}
a_n\, =\, q^{-1}\, 1
\ee
is called {\em $SL_q(n)$ type RE algebra} and denoted as ${\goth L\scriptstyle SL_q(n)}[\R]$.
\end{defin}

\begin{rem}
{\rm\small
\lb{rem-3.10}
An actual value of a numeric factor
in the right hand side of (\ref{RE-sl}) is not relevant for the definition.
Our choice
allows avoiding numeric factors later in formula
(\ref{aut1}) (see proof of proposition \ref{prop4.0.1}).

Consider realization of the RE algebra ${\goth L\scriptstyle SL_q(n)}[\R]$ as a subalgebra in
the quasi-triangular Hopf algebra ${\goth A}_{\cal R}$ (see example \ref{3.2.2}).
In this case the condition (\ref{RE-sl}) is consistent with the  pairing
$\langle\cdot ,\cdot\rangle$ of the dual Hopf algebras ${\goth A}_{\cal R}$ and ${\goth A}^*_{\cal R}$
only for the chosen normalizations (\ref{L-iz-univ}) for $L$ and $\eta=q^{1/n}$
for $\R$ (\ref{R-iz-univ}).
This point is explained in appendix \ref{append-A}, see (\ref{pair-Tan}).
}\end{rem}

\begin{rem}
{\rm\small
The $GL_q(n)$ type R-matrix $\R$ and its $SL_q(n)$ partner R-matrix
$\widetilde \R$ (\ref{partner-R}) define identical RE algebras.
}\end{rem}

In the theorem below we describe  Cayley-Hamilton and Newton identities specific to
the  $GL_q(n)$ type and Hecke type RE algebras.

\begin{theor}
Let $\R$ be a skew invertible {\rm R}-matrix of the Hecke type.
Assume that the conditions {\bf [k]} (\ref{restrict})
are fulfilled.
Then in the corresponding RE algebra ${\goth L\goth G}[\R]$
following {\em Cayley-Hamilton-Newton identities} \cite{IOPa,IOP}
\be
\lb{chn}
i_q \, \RTr{2, \dots ,i}
(A^{(i)} L_{\overline{2}} L_{\overline{3}}\dots L_{\overline{i}}) =
(-1)^{i+1} \sum_{j=0}^{i-1} (-q)^j\, a_j\, L_1^{i-j-1} \quad\forall\, 2\leq i\leq k\,
\ee
take place.
Multiplying by $L_1$ from the left and taking the {\rm R}-trace $\RTr{1}$
of these identities one obtains  {\em Newton relations}
for the sets of power sums $\{ p_i \}_{i\geq 1}$ and the set of elementary
symmetric functions
$\{ a_i \}_{i\geq 0}$ \cite{GPS}
\be
\lb{newton}
i_q  \, a_i + (-1)^i \sum_{j=0}^{i-1} \, (-q)^{j} a_j \,p_{i-j} = 0
\quad \forall\, 1\leq i\leq k\, .
\ee
Both sets $\{1, p_j \}_{j\geq 1}$ and $\{a_j\}_{j\geq 0}$ in this case  generate  the
characteristic subalgebra ${\goth C\goth h}[\R]$.

Assume additionally that $\R$ is an {\rm R}-matrix of the $GL_q(n)$ type.
Then the finite set $\{ a_i \}_{i= 0}^n $ generates the characteristic subalgebra
of the RE algebra  ${\goth L\scriptstyle GL_q(n)}[\R]$, and
following {\em Cayley-Hamilton  identity} is fulfilled \cite{GPS}:
\be
\lb{char}
\sum^{n}_{i=0}\, (-q)^i\, a_{i}\, L^{n-i}\, =\, 0 \, .
\ee
This identity leads, in particular, to an invertibility of the matrix $L$:
$$
%\be\lb{L-inv}
L^{-1}\, =\, q^{-1}\, a_n^{-1}\, \sum_{i=0}^{n-1}\, (-q)^{-i}\, a_{n-i-1}\, L^i\, .
%\ee
$$
\end{theor}

\begin{rem}
\lb{rem3.5}
{\rm\small
One can introduce generating functions $a(x)$, $p(x)$ for the elementary symmetric functions
and for the power sums
$$
a(x) := \sum_{i\geq 0} \, a_i \,x^i \, , \qquad p(x) := \sum_{i\geq 1} \, p_i \, x^i\, .
$$
The Newton relations (\ref{newton}) can be written as
a finite difference equation for the generating functions
$$
a(qx) \, p(-x)\, =\, {a(q^{-1}x) - a(q x)\over q-q^{-1}} \, .
$$
}\end{rem}

\medskip
For the $GL_q(n)$ type RE algebra we now construct its central extension by roots of the
characteristic polynomial (\ref{char}).

\begin{defin}
\lb{3.2.8}
Denote ${\goth S}_n$ a $\mathbb C$-algebra of polynomials in $n$ pairwise commuting
invertible indeterminates $\mu^{\pm 1}_\alpha$ and their differences
$(\mu_\alpha-\mu_\beta)^{\pm 1}$,
$~\alpha ,\beta =1,\dots , n$, $~\alpha\neq \beta$.

Let $\R$ be a skew invertible {\rm R}-matrix of the $GL_q(n)$ type,
${\goth L\scriptstyle GL_q(n)}[\R]$ be the corresponding RE algebra, and
${\goth C\goth h}[\R]$ be it's characteristic subalgebra.
Consider a monomorphism ${\goth C\goth h}[\R]\hookrightarrow {\goth S}_n$
defined on generators
as~\footnote{
When defining the map (\ref{Ch-E-hom}) we implicitly assume an algebraic independence
of the elements $a_i$, $i=1,\dots ,n$. Otherwise, we should impose the
same algebraic conditions on functions $e_i(\mu_1,\dots ,\mu_n)$.
}
\be
\lb{Ch-E-hom}
a_i\mapsto e_i(\mu_1,\dots ,\mu_n):= \sum_{1\leq j_1<\dots <j_i\leq n}
\mu_{j_1} \mu_{j_2}\dots \mu_{j_i}
\qquad
\forall\, i=0,1,\dots ,n\, ,
\ee
where $e_i$ are the elementary symmetric functions of their arguments.
The map (\ref{Ch-E-hom}) defines naturally a structure of, say, left ${\goth C\goth h}[\R]$-module
on  ${\goth S}_n$. A central extension of the algebra
${\goth L}{\scriptstyle GL_q(n)}[\R]$
$$
\overline{\goth L}{\scriptstyle GL_q(n)}[\R] \, :=\,
{\goth L}{\scriptstyle GL_q(n)}[\R]\,
\raisebox{-4pt}{$\bigotimes\atop {\goth C\goth h}_{\R}$}\,
{\goth S}_n\, :
$$
\be
\lb{REA-spec-ext}
a_\alpha\, =\, e_\alpha(\mu_1 ,\dots ,\mu_n)\, , \quad
L^i_j\,\mu_\alpha\, =\, \mu_\alpha\, L^i_j\qquad
\forall\,  i,j=1,\dots ,\dim V ,\;\;\forall\,
\alpha=1,\dots ,n\, ,
\ee
is called a {\em (semisimple) spectral completion of  ${\goth L}{\scriptstyle GL_q(n)}[\R$]}.
A  quotient of this algebra by  relations
$$
%\be\lb{map-sl}
a_n\, =\, \prod_{\alpha=1}^n \mu_\alpha\, =\, q^{-1}\, .
%\ee
$$
is called a {\em (semisimple) spectral completion of  ${\goth L\scriptstyle SL_q(n)}[\R]$}
 and denoted as
$\overline{\goth L}{\scriptstyle SL_q(n)}[\R]$.
Variables $\mu_\alpha$ are called {\em spectral variables}.
\end{defin}

\begin{rem}{\rm\small
Assuming that the spectral variables $\mu_\alpha$
are invariants of the coadjoint coaction,
% (\ref{REA-com}),
the algebra
$\overline{\goth L}{\scriptstyle GL_q(n)}[\R]$ ($\overline{\goth L}{\scriptstyle SL_q(n)}[\R]$)
inherits the structure of left coadjoint
${\goth F\scriptstyle GL_q(n)}[\R]$- (${\goth F\scriptstyle SL_q(n)}[\R]$-) comodule
algebra.
}
\end{rem}

\begin{cor}
\lb{cor3.16}
In the spectrally completed algebra $\overline{\goth L}{\scriptstyle GL_q(n)}[\R]$
the characteristic identity (\ref{char}) assumes a  factorized form
\be
\lb{CH-factor}
\prod_{\alpha =1}^n\bigl(L - q\mu_{\alpha} I\bigr)\, =\, 0\, .
\ee
One can construct a resolution of the matrix unity
\be
\lb{proj}
P^\alpha\, :=\, \prod_{\beta=1\atop \beta\neq\alpha}^n{\bigl
(L - q\mu_{\beta} I\bigr)\over q(\mu_\alpha-\mu_\beta)}\; :\quad
P^\alpha P^\beta\, =\, \delta_{\alpha\beta}\, P^\alpha\, , \quad
\sum_{\alpha=1}^n P^\alpha\, =\, I\, ,
\ee
\vspace{-6mm}\ni
so that
\be
\lb{L-P}
L\, P^\alpha\, =\, P^\alpha L\, =\, q\mu_\alpha P^\alpha\, .
\ee
\end{cor}

\begin{rem}{\small\rm
In papers \cite{GS1,DM1,DM2,GS2} the factorized form of
the Cayley-Hamilton identity and the
projectors $P^\alpha$ where used to construct
explicitly  quantized semisimple coadjoint orbits
of $GL(n)$ and line bundles over them.}
\end{rem}

\subsection{Quantized differential operators over matrix group
(Heisenberg double)}
\lb{subsec3.3}

\begin{defin} {\rm \cite{AF1,STS}}
\lb{3.3.1}
Let $\R$, $T$ and $L$ be as described in the definitions \ref{def3.1} and \ref{3.2.1}.
A {\em Heisenberg double (HD) algebra} ${\goth D\goth G}[\R,\gamma]$
of the two algebras ${\goth F\goth G}[\R]$ and
${\goth L\goth G}[\R]$ is an associative unital algebra
generated by the components of the matrices $T$ and $L$
subject additionally to a permutation relation
\be
\lb{LT}
\gamma^{2} \, T_1 L_2 \,=\, \R_{12} L_1 \R_{12} T_1\,   ,\quad \mbox{where~~}
\gamma\in\{{\mathbb C}\backslash 0\}\, .
\ee
The HD algebra carries structures of left and right
${\goth F\goth G}[\R]$-comodule algebra, respectively,
\ba
\lb{delta-l}
\delta_{\ell} (T^i_j)\, =\, \sum_k\, T^i_k \otimes T^k_j \, , &\quad&
\delta_{\ell} (L^i_j)\, =\, \sum_{k,m}\, T^i_k \, (T^{-1})^m_j \otimes L^k_m  \, ;
\\[2mm]
\lb{delta-r}
\delta_r (T^i_j)\, =\, \sum_k\, T^i_k \otimes T^k_j \, , &\quad&
\delta_r (L^i_j) = L^i_j \otimes 1  \, .
\ea
\end{defin}

\begin{example}{\rm %\small
\lb{example-3.19}
The Heisenberg double is closely related to a {\em smash
product} of two mutually dual Hopf algebras (see, e.g., \cite{Mon}).
Namely, given a pair
${\goth A}_{\cal R}$ and ${\goth A}_{\cal R}^*$
their smash product algebra ${\goth A}_{\cal R}\sharp {\goth A}_{\cal R}^*$
is a linear space ${\goth A}_{\cal R}\otimes {\goth A}_{\cal R}^*$
supplied with a multiplication
\be
\lb{smash}
(x\,\sharp\, u)(y\,\sharp\, v)\, :=\, \langle u_{(1)}, y_{(2)}\rangle\,
(x y_{(1)}\,\sharp\, u_{(2)} v)\, ,
\ee
where $x,y\in {\goth A}_{\cal R}$, $u,v\in {\goth A}_{\cal R}^*$, and symbols
$(x\,\sharp\, u)$, $(y\,\sharp\, v)$ denote
elements of ${\goth A}_{\cal R}\sharp {\goth A}_{\cal R}^*$.

Let us calculate in the settings of the examples \ref{3.1.2}, \ref{3.2.2}
the smash product of the elements $(T^i_j\,\sharp 1)\,=\, T^i_j$ and
$(1\,\sharp\,{L^{\scriptscriptstyle (\pm)}}^i_j)\,=\,{L^{\scriptscriptstyle (\pm)}}^i_j$
\ba
\nn
%\lb{dop2}
T_1 \, L^{\scriptscriptstyle (+)}_2 &=& L^{\scriptscriptstyle (+)}_2 \langle T_1 ,
\, L^{\scriptscriptstyle (+)}_2 \rangle \, T_1\, =\, \eta^{-1}\,
L^{\scriptscriptstyle (+)}_2 \, P_{12}{\R}_{12} \, T_1 \, ,
\\[2mm]
\nn
T_1 \, L^{\scriptscriptstyle (-)}_2 &=& L^{\scriptscriptstyle (-)}_2 \langle T_1 , \,
L^{\scriptscriptstyle (-)}_2 \rangle \, T_1\, =\, \eta\,
L^{\scriptscriptstyle (-)}_2 \, P_{12}{\R}^{-1}_{12} \, T_1 \, ,
\ea
wherefrom it follows that the smash product of $T_1$ and
$L_2\propto ({L^{\scriptscriptstyle (-)}}^{-1} L^{\scriptscriptstyle (+)})_2$
is given by  (\ref{LT}) with $\gamma =\eta$.
However, we stress that in general
one can keep $\gamma$ independent of the normalization $\eta$ of the R-matrix
at a price of loosing universality of formulas.
Indeed, the multiplication in the smash product algebra
is given by (\ref{smash}) universally for any pair of its elements, while
the relation (\ref{LT}) in the HD algebra is written for the generators $L$ and $T$ only.
}
\end{example}

Now we consider a geometric interpretation of the HD algebra.
Applying the substitution $L^i_j\rightarrow \ell^i_j$ (\ref{ell}) and taking
the "classical" limit $q\rightarrow 1$ in relation (\ref{LT}) we find
\be
\lb{ellT}
\bigl[T_1,\ell_2\bigr]\, =\, \bigl(P_{12}-\gamma'\, I_{12}\bigr)\, T_1\, ,
\ee
where
we used the Hecke condition (\ref{H}) in a form $\,\R^2 = \id + (q-q^{-1})\R\,$
and assumed additionally $\R\stackrel{q\rightarrow 1}{\longrightarrow} P\,$
(which is true for the Drinfeld-Jimbo R-matrix (\ref{DJ}))
and $\,\gamma \equiv \gamma(q) = 1 + (q-q^{-1}){\gamma'\over 2}+ o(q-q^{-1})$.
The commutation relations (\ref{ellT}), (\ref{REA-limit}) are realized by  operators
\be
\nn
%\lb{ellT-real}
\ell^i_j\, =\, \sum_{k=1}^n g^i_k {\partial\over\partial g^j_k}\, , \quad
T^i_j\, =\, |g|^{-\gamma'} g^i_j\,, \quad\mbox{where}\quad
g^i_j:={\rho_V(g)}^i_j\, ,\; g\in GL(n)\, ,\;
|g|:=\det {\cal k}g^i_j{\cal k}\, .
\ee
These are, respectively, right invariant vector fields and
properly normalized coordinate functions on $GL(n)$.
Together they generate an algebra of differential operators over $GL(n)$\footnote{
Imposing conditions ~$\gamma'=1/n$,~
~$\det T=1$, ~$\tr\ell=0$~
one can make a reduction
to a subalgebra of differential operators over $SL(n)$.}.

Extrapolating the classical picture we can treat  ${\goth D\goth G}[\R,\gamma]$
as an algebra of {\em quantized differential operators over matrix group} or,
equivalently, as {\em quantized functions over cotangent bundle of a matrix group}
(see \cite{AF1,AF2,SWZ1,IP}).
The form of the substitution (\ref{ell}) suggests that the quantized vector fields $L^i_j$
possess properties of finite difference operators rather then of the differential operators.
In particular, they do not satisfy classical Leibniz rule when acting on
functions (see (\ref{LT})).
\medskip

The next proposition describes an action of the characteristic
subalgebra on quantized functions in the Hecke case.

\begin{prop}
Let $\R$ be a skew invertible Hecke type
{\rm R}-matrix.
Assume that  the conditions {\rm\bf  [k]}
(\ref{restrict}) are satisfied, so that
the elementary symmetric functions
$a_i\in {\goth C\goth h}[\R]
%\subset {\goth L\goth G}[\R]
\subset {\goth D\goth G}[\R,\gamma]$,
$0\leq i\leq k$,
(\ref{sigma}) are well defined.
Then relations
\be
\lb{Tsigma}
\gamma^{2i} T \, a_i \, =\, a_i \, T\, -\,
(q^2-1) \sum_{j=1}^{i} (-q)^{-j} a_{i-j}\, (L^j  T) \qquad\forall\, 0\leq i\leq k
\ee
are fulfilled for the {\em Hecke type HD
algebra} ${\goth D\goth G}[\R,\gamma]$.
\end{prop}

\ni
{\em Proof.~}
For any operator $Y\in {\rm End}_W(V^{\otimes i})$, where $W$ is an arbitrary
$\mathbb C$-linear space, we denote
\be
\lb{up-1}
Y^{\uparrow 1}\, :=\, (P_1 P_2\dots P_i) Y  (P_1 P_2\dots P_i)^{-1}\, .
\ee
For any R-matrix $\R$ we define series of operators $J_i$,~
$Z_i$
\be
\lb{J-k}
J_1 := I\, ,\qquad J_{i+1}:= \R_i\, J_i\, \R_i\qquad \forall\, i\geq 1\, ,
\qquad Z_i:= \prod_{j=1}^i J_j\, .
\ee
\vspace{-5mm}
\begin{rem}
\lb{rem3.3.4}
{\small\rm
Elements $J_j$, $1\leq j\leq i$, are R-matrix realizations of a remarkable
set of {\em Jucys-Murphy
elements} in the braid group ${\cal B}_i$:
\be
\nn
%\lb{JM}
j_1:=1\, , \qquad j_{j+1}:= \sigma_j\, j_j\, \sigma_j\qquad\forall\, j=1,\dots ,i-1\, .
\ee
These elements generate a commutative subgroup in ${\cal B}_i$ and their product
$z_i:=\prod_{j=1}^i j_j$ is a central element in  ${\cal B}_i$.
For their applications and for historical references see, e.g., \cite{OP-Lect}.
}
\end{rem}

With these notations permutation relations (\ref{REA}), (\ref{LT}) can be suitably written
for arbitrary R-copies of the matrix $L$:
\ba
\lb{LiLj}
(L_{\overline{i}}J_i)(L_{\overline{j}}J_j)&=&(L_{\overline{j}}J_j)(L_{\overline{i}}J_i),
\\[1mm]
\lb{LkT}
\gamma^2\, T_1 (L_{\overline i} J_i)^{\uparrow 1}& =&
(L_{\overline{i+1}}J_{i+1})\, T_1\, =\, \R_i\, (L_{\overline{i}}J_i)\,\R_i\, T_1\qquad
\forall\, i,j\geq 1\, .
\ea
Here the second equality follows from the recursive definitions of $L_{\overline{i+1}}$
and $J_{i+1}$, while the first equality can be easily proved by induction on $i$.

Next, we prepare a suitable  expression  for $a_i$ (\ref{sigma}):
\ba
\nn
a_i &=& \RTr{1,\dots ,i}\Bigl(L_{\overline{1}}\dots L_{\overline{i}}\, A^{(i)}\Bigr)\; =\;
q^{i(i-1)}\,\RTr{1,\dots ,i}\Bigl(L_{\overline{1}}\dots L_{\overline{i}}\, Z_i\, A^{(i)}\Bigr)
\\[2mm]
\lb{a-alt}
&=&
q^{i(i-1)}\,\RTr{1,\dots ,i}\Bigl((L_{\overline{1}}J_1)\dots
(L_{\overline{i}} J_i)\,  A^{(i)}\Bigr)
\ea
Here we substituted $Z_i A^{(i)} = q^{-i(i-1)} A^{(i)}$ in the first line
and used a commutativity
relation
\be
\lb{JL}
L_{\overline{i}}\, J_j\, =\, J_j\, L_{\overline{i}}\qquad\forall\,i,j:\,  i>j\,
\ee
in the second line.
By relabelling  the subscript indices of the R-traces we then recast
(\ref{a-alt}) in a following form\footnote{
Notice a similarity of the formula (\ref{a-altern}) with the relation (\ref{dmat4}).
The role of the R-matrices $R^{\pm \varepsilon}$ is now played by the permutation
matrix $P$ (see (\ref{up-1})).}
\be
\lb{a-altern}
I_1\, a_i
\,=\,
q^{i(i-1)}\,\RTr{2,\dots ,i+1}\Bigl((L_{\overline{1}}J_1)\dots
(L_{\overline{i}} J_i)\,  A^{(i)}\Bigr)^{\uparrow 1}\, .
\ee

Now we are ready to permute $T_1$ and $a_i$. Substituting expression (\ref{a-altern})
for $a_i$ and using relations (\ref{LkT}) and (\ref{JL}) we calculate
\ba
\nn
\gamma^{2i} T_1\, a_i&=&
\gamma^{2i} T_1\, ( I_1\, a_i)
\, =\,
q^{i(i-1)}\gamma^{2i}\, \RTr{2,\dots ,i+1} T_1 \Bigl(
(L_{\overline{1}}J_1)\dots
(L_{\overline{i}} J_i)\,  A^{(i)}\Bigr)^{\uparrow 1}
\\[2mm]
\nn
&=&
q^{i(i-1)} \RTr{2,\dots ,i+1}\Bigl((L_{\overline{2}} J_2)\dots (L_{\overline{i+1}} J_{i+1})
\, A^{(i)\uparrow 1}\Bigr)\, T_1
\\[2mm]
\nn
&=&
q^{i(i-1)} \RTr{2,\dots ,i+1}\Bigl((L_{\overline{2}}\dots L_{\overline{i+1}})\, Z_{i+1}
\, A^{(i)\uparrow 1}\Bigr)\, T_1\, .
\ea
To continue the calculation we need following formula
$$
Z_{i+1}\, A^{(i)\uparrow 1}\, =\, A^{(i)\uparrow 1}\,
Z_{i+1}\, A^{(i)\uparrow 1}\, =\, q^{-i(i-1)}\left( q^2 A^{(i)\uparrow 1}\, -\,
q^{-i}(q^2-1)  (i+1)_q\, A^{(i+1)}\right)\, ,
$$
which follows by a combination of the definitions (\ref{q-anti2}), (\ref{A-k}), (\ref{J-k}),
and relations (\ref{idemp-1}),
(\ref{idemp-2}),  (\ref{H}).
So we finish the calculation
\ba
\nn
\gamma^{2i} T_1\, a_i &=&
\RTr{2,\dots ,i+1}\Bigl[(L_{\overline{2}}\dots L_{\overline{i+1}})\,
\Bigl(q^2 A^{(i)\uparrow 1}\, -\,
q^{-i}(q^2-1)  (i+1)_q\, A^{(i+1)}\Bigr)\Bigr]\, T_1
\\[2mm]
\nn
&=&
q^2 a_i\, T_1 \, +\, (-q)^{-i}(q^2-1)\,\sum_{j=0}^i (-q)^{-j}\, a_j \bigl( L^{i-j}\, T\bigr)_1
\\[0mm]
\lb{ahaha}
&=&
a_i\, T_1\, -\, (q^2-1)\, \sum_{j=1}^i (-q)^{-j} a_{i-j} ( L^j T)_1\, .
\ea
Here we calculate the first summand in the second line taking into account equality
$$
(L_{\overline{2}}\dots L_{\overline{i+1}})\,A^{(i)\uparrow 1}\, =\,
(\R_1\dots \R_i)(L_{\overline{1}}\dots L_{\overline{i}})\,A^{(i)}
(\R_1\dots \R_i)^{-1}
$$
and using $i$ times formula (\ref{dmat4}). For calculation of the second summand we use
the Cayley-Hamilton-Newton identity (\ref{chn}). Thus (\ref{Tsigma}) is proved.
\hfill$\blacksquare$

\begin{rem}{\rm\small
For the set of power sums (\ref{s}) the permutation relations with $T^i_j$
in the Hecke case read
\be
\nn
%\lb{T-power}
\gamma^{2i} T \, p_i\, =\, p_i \, T\, +\,
(q-q^{-1})^2 \sum_{j=1}^{i-1} {(2j)_q \over 2_q } p_{i-j} \, (L^j  T)\, +\,
(q-q^{-1}) { (2i)_q \over 2_q } (L^i T) \, .
\ee
One can derive this formula  applying  the R-trace $\RTr{2}$
to an equality
$
\gamma^{2i} T_1 (L_2)^i\, =\, (\R L_1 \R)^i T_1
$
and taking into account relations
\vspace{-3mm}
\ba
\nn
(\R L_1 \R)^i& =& \R (L_1)^i \R\, +\, (q-q^{-1})
\sum_{j=1}^{i-1} \R^{2j} (L_1)^{i-j} \R (L_1)^j \, ,
\\[1mm]
\nn
\R^{2j} &=& 2_q^{-1}\left((q^{2j-1} + q^{-2j+1}) \id +
(q^{2j} - q^{-2j}) \, \R\right)  \, ,
\ea
These relations, in turn, follow inductively from the Hecke condition (\ref{H}) and the
reflection equation (\ref{REA}).
Note that in this case there is no need to impose restrictions (\ref{restrict})
on $q$.
}\end{rem}

\begin{prop}
Let $\R$ be a skew invertible $GL_q(n)$ type {\rm R}-matrix.
An extension of the corresponding HD algebra ${\goth D\goth G}[\R,\gamma]$
by the elements $(\Rdet T)^{-1}$ and $(a_n)^{-1}$, satisfying
relations
\ba
\lb{detT-L}
\gamma^{2n}\,  L\, (\Rdet T)^{-1} &=&
q^2 (\Rdet T)^{-1}\, (\O L\O^{\!\!\!-1})\, ,
\\[2mm]
\lb{T-an}
\gamma^{2n}\, (a_n)^{-1}\, T &=& q^2 T\, (a_n)^{-1}\, ,
\ea
in addition to those, given in definitions \ref{3.1.3-cor} and \ref{3.2.4},
is called {\em $GL_q(n)$ type HD algebra} and denoted as
${\goth D}{\scriptstyle GL_q(n)}[\R,\gamma]$.
% Relations (\ref{Tsigma}) in this case can be equivalently written as
% \be
% \lb{Tsigma1}
%  \gamma^{2i} T \, a_i\, =\, q^2  a_i T\, +\, (q^2-1) \sum_{j=1}^{n-i} \, (-q)^{j} \,
%  a_{i+j} (L^{-j}  T) \qquad\forall\,
% 0\leq i\leq p\, .
% \ee

Let $\R$ be a skew invertible $SL_q(n)$ type {\rm R}-matrix.
In the corresponding HD algebra ${\goth D\goth G}[\R,\gamma]$ let us
restrict the parameters by condition $\gamma^{n}=q$ and take
a quotient by relations $\Rdet T =1$ and $a_n=q^{-1} 1$.
The quotient algebra is called {\em $SL_q(n)$ type HD algebra} and
denoted as ${\goth D}{\scriptstyle SL_q(n)}[\R]$.
\end{prop}

\begin{rem}{\rm \small
Notice  consistency of the $SL_q(n)$ reduction condition $\gamma^{n} =q$
with the parameter restrictions $\eta =q^{1/n}$ in the example \ref{2.4.4}
and $\gamma =\eta$ in the example \ref{example-3.19}.
}
\end{rem}

\ni
{\em Proof.~}
Relations (\ref{detT-L}) and (\ref{T-an}) should be consistent
with permutation relations for $\Rdet T$ and $a_n$ in the algebra
${\goth D}{\goth G}[\R,\gamma]$.
Permutation relation for $a_n$ with $T$ were in fact derived
in the first line of the calculation (\ref{ahaha}) (put $i=n$ and
take into account that $A^{(n+1)}=0$ in the $GL_q(n)$ case).
Permutation relation for  $\Rdet T$ with $L$ can be derived by the same method as
for $\Rdet T$ with $T$ (see \cite{Gur}, sec.5,
or \cite {I}, calculation (3.5.39)). Given these results the consistency is obvious.

%The relation (\ref{Tsigma1}) follows by a combination of (\ref{char}) and (\ref{Tsigma}).

In the $SL_q(n)$ case ($\O\propto I$, $\gamma^{n}=q$)
the elements $\Rdet T$ and $a_n$ are central. Hence,
${\goth D}{\scriptstyle SL_q(n)}[\R]$ is consistently defined.
\hfill$\blacksquare$

\begin{cor}
In the $GL_q(n)$ type HD algebra  elements
of the characteristic subalgebra  satisfy following commutation relations with
$\Rdet T$
\be
\nn
%\lb{detT-char}
\gamma^{2nk}\, \Rdet T\, ch(x^{(k)})\, =\, q^{2k}\, ch(x^{(k)})\, \Rdet T
\quad \forall \, x^{(k)}\in {\cal H}_k(q), \; k=1,2,\dots\, .
\ee
\end{cor}

\ni
{\em Proof.~} A proof is a direct calculation of permutation
of $ch(x^{(k)})$ (\ref{center}) and $\Rdet T$ exploiting relations
(\ref{detT-L}) and properties of the matrix $\O$ (\ref{O})
\be
\nn
%\lb{prop-O}
\R_{12}\, \O_1\, \O_2\, =\, \O_1\, \O_2\, \R_{12}\, , \qquad
\O\, \D\, =\, \D\, \O\,
\ee
The latter relations are proved in \cite{OP}, section 5.3.
\hfill$\blacksquare$

\begin{theor}
\lb{main1}
Let $\R$ be a skew invertible $GL_q(n)$ ($SL_q(n)$) type {\rm R}-matrix. An extension of the
corresponding HD algebra ${\goth D}{\scriptstyle GL_q(n)}[\R,\gamma]$
(${\goth D}{\scriptstyle SL_q(n)}[\R]$) by the algebra ${\goth S}_n$
of polynomials in mutually commuting indeterminates
$\mu_\alpha^{\pm 1}$, $(\mu_\alpha -\mu_\beta)^{\pm 1}$
satisfying relations (\ref{REA-spec-ext}) together with
\be
\lb{T-mu}
\gamma^2 \, (P^\beta T)\, \mu_\alpha\, =\,
q^{2\delta_{\alpha\beta}}\mu_\alpha\, (P^\beta T)\qquad \forall\,
\alpha ,\beta =1,\dots ,n\, ,
\ee
or, equivalently,
\vspace{-2mm}
\be
\nn
%\lb{T-mu2}
\gamma^2 \,  T\, \mu_\alpha\, =\,
\mu_\alpha\, T\, +\, (q^2-1)\mu_\alpha\, (P^\alpha T)\, ,
\ee
is called a {\em (semisimple) spectral completion} of the $GL_q(n)$ ($SL_q(n)$) type HD algebra
and denoted as
$\overline{{\goth D}}{\scriptstyle GL_q(n)}[\R,\gamma]$
($\overline{{\goth D}}{\scriptstyle SL_q(n)}[\R]$).
\end{theor}

\begin{rem}{\small\rm
To avoid problems with permutations of  $(\mu_\alpha -\mu_\beta)^{-1}$
with $P^\sigma T$ one could assume invertibility of all elements
$(\mu_\alpha - q^{2k}\mu_\beta)\;\,\forall\,\alpha\neq\beta ,\, k\in\mathbb{Z}$.
Further on we will not make such permutations and so we don't impose the
corresponding restrictions.
}
\end{rem}

\begin{rem}{\small\rm
Assuming that the spectral variables $\mu_\alpha$
are invariants of
both left and right coactions,
the algebra
$\overline{\goth L}{\scriptstyle GL_q(n)}[\R,\gamma]$ ($\overline{\goth L}{\scriptstyle SL_q(n)}[\R]$)
inherits the structures of left and right
${\goth F\scriptstyle GL_q(n)}[\R]$- (${\goth F\scriptstyle SL_q(n)}[\R]$-) comodule
algebra (see definition \ref{3.3.1}).
}
\end{rem}

\begin{rem}{\small\rm
Note that relation (\ref{T-mu}) is typical for Weyl algebra generators.
In fact there are many ways to combine from the elements $(P^\beta T)_{ij}$
a set of $n$ generators satisfying  Weyl relations
with the spectral variables $\mu_\alpha$.
One such possibility is used later in section \ref{subsec4.4}.
}
\end{rem}

\ni
{\em Proof.~}
We have to check consistency of relations (\ref{T-mu}), (\ref{Tsigma}) with the conditions
$a_i=e_i(\mu_1,\dots ,\mu_n)\equiv e_i(\mu)$ for $1\leq i\leq n$.
Denote
$e_i(\mu^{\,/ \alpha}):=e_i(\mu)|_{\mu_\alpha=0}\,$. We have
\be
\lb{mu-strih}
e_i(\mu)\, =\, e_i(\mu^{\, /\alpha})\, +\, \mu_\alpha\, e_{i-1}(\mu^{\, /\alpha})\quad
\Rightarrow\quad
e_i(\mu^{\, /\alpha})\, =\, \sum_{j=0}^i (-\mu_\alpha)^j\, e_{i-j}(\mu)\, .
\ee
Using relations (\ref{T-mu}), (\ref{mu-strih}), (\ref{proj}) and (\ref{L-P}) we calculate
\ba
\nn
\gamma^{2i}\, T\, e_i(\mu) &=&
\gamma^{2i}\sum_{\alpha=1}^n (P^\alpha T)\left(e_i(\mu^{\, /\alpha})
\, +\, \mu_\alpha\, e_{i-1}(\mu^{\, /\alpha})\right)
\\
\nn
&=& \sum_{\alpha=1}^n\left(e_i(\mu^{\, /\alpha})
\, +\, q^2\mu_\alpha\, e_{i-1}(\mu^{\, /\alpha})\right) (P^\alpha T)
\\
\nn
&=& \sum_{\alpha=1}^n\Bigl(e_i(\mu)\, +\, (q^2-1)\,\mu_\alpha \sum_{j=0}^{i-1}
(-\mu_\alpha)^j e_{i-j-1}(\mu)\Bigr) (P^\alpha T)
\\
\nn
&=& \Bigl(e_i(\mu)\, -\, (q^2-1)\sum_{j=1}^i (-L/q)^j e_{i-j}(\mu)\Bigr)
\sum_{\alpha=1}^n (P^\alpha T)
\\
\nn
&=& e_i(\mu)\,T\, -\, (q^2-1)\sum_{j=1}^i (-q)^{-j} e_{i-j}(\mu)\, (L^j T)\, ,
\ea
which coincides with (\ref{Tsigma}) under identification $e_i(\mu)=a_i$.
\hfill$\blacksquare$

\begin{cor}
\lb{cor3.25}
In the completed $GL_q(n)$ type HD algebra $\overline{{\goth D}}{\scriptstyle GL_q(n)}[\R,\gamma]$
following permutation relations hold
\be
\lb{detT-mu}
\gamma^{2n}\, \Rdet T\, \mu_\alpha\, =\, q^{2}\, \mu_\alpha\, \Rdet T
\quad \forall \, \; \alpha=1,2,\dots ,n\, .
\ee
\end{cor}

\ni
{\em Proof.~}
Using formulas (\ref{detT2}), (\ref{proj}), (\ref{T-mu}) we
can permute $\Rdet T$ and $\mu_\alpha$:
\ba
\nn
\gamma^{2n} \Rdet T \mu_\alpha &\!\!\!\!=&\!\!\!\!
\gamma^{2n}\!\!\!\!
\sum_{\beta_1,\dots ,\beta_n=1}^n \Tr{1,\dots ,n} \bigl(A^{(n)}
(P^{\beta_1} T)_1\dots (P^{\beta_n} T)_n\bigr) \mu_\alpha
\\[1mm]
\lb{detT-mu-promezh}
&\!\!\!\!=&\!\!\!\!
\mu_\alpha\!\!\!\!
\sum_{\beta_1,\dots ,\beta_n=1}^n q^{2\sum_{j=1}^n\delta_{\alpha\beta_j}}
\Tr{1,\dots ,n} \bigl(A^{(n)}
(P^{\beta_1} T)_1\dots (P^{\beta_n} T)_n\bigr) .
\ea
Assuming that
\be
\lb{proj-zero}
\Tr{1,\dots ,n} \bigl(A^{(n)}
(P^{\beta_1} T)_1\dots (P^{\beta_n} T)_n\bigr)\, =\, 0\, ,
\quad \mbox{if there exists a pair}\; i,j:\; \beta_i=\beta_j\, ,
\ee
we conclude that for any nonzero summand in (\ref{detT-mu-promezh})
the coefficient $q^{2\sum_{j=1}^n\delta_{\alpha\beta_j}}$ equals $q^2$, and therefore
we can complete the calculation
$$
\gamma^{2n}\, \Rdet T\, \mu_\alpha \, =\,
q^2 \mu_\alpha\,\sum_{\beta_1,\dots ,\beta_n=1}^n
\Tr{1,\dots ,n} \bigl(A^{(n)}
(P^{\beta_1} T)_1\dots (P^{\beta_n} T)_n\bigr)\, =\,
q^2\, \mu_\alpha\, \Rdet T\, .
$$

It lasts to prove the assumption. First, we note that conditions on $\beta_i$ in (\ref{proj-zero})
stand that there exists integer $\sigma$: $1\leq \sigma\leq n$, and $\sigma\neq\beta_i\;\forall\,i$.
Therefore, any projector $P^{\beta_i}$ in (\ref{proj-zero}) contains factor $(L-q\mu_\sigma I)$.
Using relations (\ref{LkT}), (\ref{T-mu}) we can move all such factors to the left side
of the expression. Thus we obtain
\ba
\lb{lhs_of}
\mbox{left hand side of (\ref{proj-zero})}\, \propto
\Tr{1,\dots ,n} \left(A^{(n)}\left\{{\textstyle \prod_{j=1}^n}\bigl(L_{\overline{j}}J_j-q\mu_\sigma I)
\right\}\dots \right) .
\ea
Next, we note
that  the expression in braces is a symmetric function in a commuting set of matrices
$L_{\overline{j}}J_j$ (see (\ref{LiLj}))
which by relations
$$
R_i (L_{\overline{i}}J_i)(L_{\overline{i+1}}J_{i+1}) =
(L_{\overline{i}}J_i)(L_{\overline{i+1}}J_{i+1}) R_i ,\quad
R_i (L_{\overline{i}}J_i + L_{\overline{i+1}}J_{i+1}) =
(L_{\overline{i}}J_i + L_{\overline{i+1}}J_{i+1}) R_i ,
$$
and by
(\ref{LR-comm})
together with the same formulas for $J_k$
commutes with  $R_i,\; i=1,\dots ,n-1$, and so
with $A^{(n)}$.
Hence, using relations $A^{(n)}=(A^{(n)})^2$ and $\mbox{rk} A^{(n)}=1$
we can separate a left factor
$
\kappa :=\Tr{1,\dots ,n} \bigl(A^{(n)}\prod_{j=1}^n\bigl(L_{\overline{j}}J_j-q\mu_\sigma I)
\bigr)
$
in (\ref{lhs_of}). This factor we now calculate explicitly.

Taking into account relations (\ref{D-A}), (\ref{LiLj}), (\ref{JL}) and $A^{(n)}J_i=q^{-2(i-1)}A^{(n)}$
we transform the expression for $\kappa$:
\be
\lb{kappa}
\kappa\, =\, q^n\,
\RTr{1,\dots ,n} \left(A^{(n)}{\textstyle \prod_{j=1}^n}\bigl(L_{\overline{j}}\,-\,q^{2j-1}\mu_\sigma I)
\right).
\ee
Expanding this expression in powers of $L$ and noticing that
(\ref{A-A}) assumes $\RTr{k+1,\dots ,n} A^{(n)}= q^{n(k-n)}\frac{(k)_q!(n-k)_q!}{(n)_q!}:=
q^{n(k-n)}{n\choose k}^{-1}$
we find that $k$-th order monomials
$$
\RTr{1,\dots ,n} \left(A^{(n)}L_{\overline{i_1}}\dots L_{\overline{i_{k}}}\right)
= \RTr{1,\dots ,n} \left(A^{(n)}L_{\overline{1}}\dots L_{\overline{k}}\right)\, =\,
q^{n(k-n)}{\textstyle {n\choose k}^{-1}_q}\, a_k\,
$$
are equal to each other for any choice of indices $1\leq i_1<\dots i_k\leq n$.
Their corresponding coefficients in (\ref{kappa}) sum up to
$$
(-q^{-1}\mu_\sigma)^{n-k}
\sum_{1\leq i_1<\dots i_{n-k}\leq n} q^{2\sum_{r=1}^{n-k} i_r }\, =\, q^{n(n-k)}
{\textstyle {n \choose k}}_q (-\mu_\sigma)^{n-k}\, ,
$$
and so we obtain
$$
\kappa\, =\, q^n\, \sum_{k=0}^n a_k (-\mu_\sigma)^{n-k}\, =\,
q^n\, \prod_{\alpha=1}^n(\mu_\alpha -\mu_\sigma)\, =\, 0\, ,
$$
where we took (\ref{REA-spec-ext}) into account.
\hfill$\blacksquare$

\subsection{Quantized left invariant vector fields.}
\lb{subsec3.4}

In a classical differential geometry of the Lie groups one uses
two global bases on tangent bundles -- the bases of right and left invariant vector fields.
In previous sections we discussed quantization of the right invariant vector fields only
and defined the HD algebra ${\goth D\goth G}[\R,\gamma]$ in their terms.
To demonstrate a left-right symmetry of the whole construction we now describe the HD algebra
using a set of left invariant generators. We also find explicit relations between
the spectra of left and right invariant vector fields in both, the $GL_q(n)$
and the $SL_q(n)$ cases.

In the assumptions of the definition \ref{3.3.1} consider a matrix
$M$ whose components belong to
${\goth D\goth G}[\R,\gamma]$:
\be
\lb{M}
M^i_j\, :=\, \sum_{k,m}\,(T^{-1})^i_k\, L^k_m\, T^m_j\, .
\ee
Taking into account transformation properties of the matrix elements $M^i_j$
with respect to the left and right ${\goth F\goth G}[\R]$-coactions (\ref{delta-l})
and (\ref{delta-r})
\be
\lb{delta-M}
\delta_\ell (M^i_j)\, =\, 1\otimes M^i_j\, , \qquad
\delta_r(M^i_j)\, =\, \sum_{k,m}\, (T^{-1})^i_k\, T^m_j \otimes M^k_m\, ,
\ee
we shall call them a {\em basis of  quantized left invariant vector fields over matrix group}.

One can give presentation of the HD algebra ${\goth D\goth G}[\R,\gamma]$ in terms of
generators~ $T^i_j$ and $M^i_j$,~ and relations
\ba
\nn
%\lb{Rdot-TT}
\Rd_{12}\, T_2\, T_1  &=& T_2\, T_1 \Rd_{12}\, ,
\\[0mm]
\lb{REA-M}
\Rd_{12}\, M_1\,\Rd_{12}\, M_1  &=& M_1\,\Rd_{12}\,M_1\,\Rd_{12}\, ,
\\[0mm]
\lb{RMRT}
\gamma^{-2}\, M_2\, T_1& =& T_1\, \Rd_{12}\,M_1\,\Rd_{12}\, ,
\ea
where we denote
\be
\lb{PRP}
\Rd_{12}\, :=\, (P\R^{-1}P)_{12}\, =\, (\R_{21})^{-1}\, .
\ee
Necessary technical data about $\Rd$ are collected in  lemma \ref{lem-a1} in  appendix
\ref{append}.
\smallskip

By (\ref{REA-M}), the entries of matrix $M$  generate yet another
RE subalgebra ${\goth L\goth G}[\Rd]$\vspace{-0.9mm} in the HD algebra ${\goth D\goth G}[\R,\gamma]$.
By (\ref{delta-M}), the subalgebra ${\goth L\goth G}[\Rd]$ is
a right coadjoint ${\goth F\goth G}[\R]$-comodule algebra.
We also notice a nontrivial but quite expected property of the quantized left and right
invariant vector fields --- their mutual commutativity
\be
\nn
%\lb{LM}
M_1\, L_2\, =\, L_2\, M_1\, .
\ee

In the rest of this section we
investigate the characteristic subalgebra\vspace{-0.7mm}
${\goth C\goth h}[\Rd]\subset {\goth L\goth G}[\Rd]$. In particular, we shall see that
${\goth C\goth h}[\Rd]={\goth C\goth h}[\R]$ for the ${\goth D\goth G}[\R,\gamma]$-subalgebras
${\goth L\goth G}[\Rd]$ and
${\goth L\goth G}[\R]$.%\vspace{-0.7mm}

It is suitable to introduce {\em $\Rd$-copies}
of the matrix $M$ (c.f. with (\ref{L-over}))
\be
\lb{M-copies}
\Md{1}\, :=\, M_1\, , \qquad \Md{k+1}\, :=\, \Rd_k\, \Md{k}\, (\Rd_k)^{-1}\, ,\vspace{-3mm}
\ee
and $\Rd$-matrix realizations
of the Jucys-Murphy elements (c.f. with (\ref{J-k}) and remark \ref{rem3.3.4})
\be
\lb{Jd}
\Jd_{1}\, :=\, I\, , \qquad \Jd_{k+1}\, :=\, \Rd_k \Jd_k \Rd_k
\qquad \forall\, k\geq 1\, .
\ee
In their terms the relations (\ref{REA-M}), (\ref{RMRT})
can be  written as (c.f. with (\ref{REA-forms}), (\ref{LkT}))
\ba
\nn
%\lb{RMM-copies}
\Rd\,\Md{k}\,\Md{k+1}
& =&
\Md{k}\,\Md{k+1}\,\Rd \, ,
\\[0mm]
\lb{MT-copies}
\gamma^{-2}\, \Bigl(\Md{k}\,\Jd_k\Bigr)^{\uparrow 1}\, T_1& =&
T_1 \Bigl(\Md{k+1}\,\Jd_{k+1}\Bigr)\, .
\ea

We now assume that the R-matrix $\Rd$ is skew invertible\footnote{This is indeed the case if
$\R$ is skew invertible and $\D$ is invertible (see lemma \ref{lem-a1} in appendix \ref{append}).}
and introduce two
generating sets in the characteristic subalgebra
 ${\goth C\goth h}[\Rd]\subset {\goth L\goth G}[\Rd]$:
the power sums
$\stackrel{\ast}{p}_k$
\be
\nn
%\lb{n-d}
\stackrel{\ast}{p}_k \, :=\, \Rdtr (M^k)\, ,\qquad k=1,2,\dots\, ,
\ee
and, assuming additionally that conditions {\rm\bf  [k]} (\ref{restrict})
are fulfilled, the elementary symmetric functions $\stackrel{\ast}{a}_i$
\be
\lb{a-d}
\stackrel{\ast}{a}_0 \, :=\, 1\, , \qquad
\stackrel{\ast}{a}_i \, :=\, \RdTr{1,\dots ,i}
\Bigl(\Adro{(i)} \Md{1}\Md{2}\dots\Md{i} \Bigr)\qquad \forall\, 1\leq i\leq k\, .
\ee

\begin{prop}
\lb{prop3.4.2}
Let $\R$ be a skew invertible Hecke type
{\rm R}-matrix and $\D$ be invertible.
Assume that conditions {\rm\bf  [k]}
(\ref{restrict}) are satisfied.
Then for two sets of elements in
${\goth C\goth h}[\R]\subset {\goth D\goth G}[\R,\gamma]$ ---
$a_i$ (\ref{sigma}) and $\stackrel{\ast}{a}_i$ (\ref{a-d}) ---
following relations are satisfied
\be
\nn
%\lb{a-a-ast}
\stackrel{\ast}{a}_i\, =\, \gamma^{2i}\, a_i\qquad \forall\, 0\leq i\leq k\, .
\ee
\end{prop}

\ni
{\em Proof.~} We transform the expression (\ref{a-d}) for $\stackrel{\ast}{a}_i$
in a following way
\ba
\nn
\stackrel{\ast}{a}_i\, =\, \RdTr{1,\dots ,i}
\Bigl( \Md{1}\dots\Md{i} \Adro{(i)}\Bigr)\, =\,
q^{-i(i-1)}\,\RdTr{1,\dots ,i}
\Bigl( (\Md{1}\Jd_1)\dots(\Md{i}\Jd_i) \Adro{(i)}\Bigr)
\ea
Here we used formulas
\be
\lb{tech-vspom}
\Jd_k\!\! \Adro{(i)} = q^{2(k-1)}\!\!\!\Adro{(i)}\;\forall\,1\leq k\leq i\, ,\quad \mbox{~and}\quad
\Md{i}\Jd_k=\Jd_k\Md{i}\;\forall\,1\leq k< i\, ,
\ee\vspace{-5mm}

\ni
which, in turn, follow from  (\ref{R-Adro}), (\ref{M-copies}), (\ref{Jd}).

Then we apply lemma \ref{lem-a2} from the appendix \ref{append}
and use the relations (\ref{P-Theta}) to move
$\Adro{(i)}$ leftwards
\ba
\nn
&&
\hspace{-6mm}\stackrel{\ast}{a}_i\, =
q^{-i(i-1)}\,\gamma^{i(i+1)}\,\RdTr{1,\dots ,i}\RTr{i+1,\dots, 2i}
\Bigl(  \Upsilon^{(i)}_P \Upsilon^{(2i)}_P (LT)_{i}\dots
(LT)_1 \Upsilon_{\subRd}^{(2i)}\!\!\Adro{(i)\uparrow i}\times
\\
\nn
&&\hspace{100mm}
\times (T_{i}\dots T_1)^{-1}
\Upsilon^{(2i)}_P \Upsilon^{(i)}_P\Bigr),
\ea\vspace{-3mm}

\ni
Here matrices $\Upsilon_*^{(*)}$ are defined in (\ref{Theta}).

Next, we  permute $\Adro{(i)\uparrow i}$ with $\Upsilon_{\subRd}^{(2i)}$ and
cancel terms $\Upsilon^{(i)}_P \Upsilon^{(2i)}_P$ on the left and
$\Upsilon^{(2i)}_P \Upsilon^{(i)}_P$ on the right. The latter cancellation exchange
the R-traces $\Rtr$ and $\Rdtr$:
\ba
\nn
%\hspace{-12mm}
\stackrel{\ast}{a}_i\,=\,
q^{-i(i-1)}\,\gamma^{i(i+1)}\,\RTr{1,\dots ,i}\RdTr{i+1,\dots, 2i}
\Bigl((LT)_{i}\dots(LT)_1 \Adro{(i)}\!\Upsilon_{\subRd}^{(2i)}
(T_{i}\dots T_1)^{-1}\Bigr).&&
\ea
In the resulting expression all the R-traces $\RdTr{i+1,\dots ,2i}$ can be evaluated
with the help of lemma \ref{lem-a3}. So, we continue
\ba
\nn
&&
\hspace{-12mm}\stackrel{\ast}{a}_i\,=
q^{i(i-1)}\,\gamma^{i(i+1)}\,\RTr{1,\dots ,i}
\Bigl((LT)_{i}\dots(LT)_1 \Adro{(i)}
(T_{i}\dots T_1)^{-1}\Bigr)
\\[2mm]
\nn
&&
\hspace{-8mm}=
q^{i(i-1)}\,\gamma^{i(i+1)}\,\RTr{1,\dots ,i}
\Bigl((LT)_{1}\dots(LT)_i(T_{1}\dots T_i)^{-1}A^{(i)}\Bigr).
\ea\vspace{-5mm}

\ni
Here in the first line we used formula (\ref{Theta-J}) and (\ref{tech-vspom});
in the second line we applied formula (\ref{RdroA}) and then,
moved the two terms $\Upsilon^{(i)}_P$,
respectively,
to the left and to the right
and cancelled them under the R-traces $\RTr{1,\dots ,i}$.
Finally, using repeatedly permutation relations (\ref{LkT}) and then
formula (\ref{a-alt})
we complete the transformation
\ba
\nn
&&
\stackrel{\ast}{a}_i\,=
q^{i(i-1)}\,\gamma^{i(i+1)}\,\RTr{1,\dots ,i}
\Bigl((LT)_{1}\dots(LT)_{i-1} \Bigl((L_{\overline{1}}J_1)^{\uparrow 1} T_1^{-1}\Bigr)^{\uparrow (i-2)}
(T_{1}\dots T_{i-2})^{-1}A^{(i)}\Bigr)\hspace{1mm}
\\[1mm]
\nn
&&
\hspace{3mm}\dots =
q^{i(i-1)}\,\gamma^{2i}\,\RTr{1,\dots ,i}
\Bigl((L_{\overline{1}}J_1)\dots(L_{\overline{i}}J_i)
A^{(i)}\Bigr)\, =\,
\gamma^{2i}\, a_i\, .\hspace{44mm}\blacksquare
\ea

\begin{rem}{\rm\small
For the sets of power sums $p_i$ and \raisebox{-1pt}{$\stackrel{\ast}{p}_i$} one can prove
following recurrent relations
\be
\nn
%\lb{n-p}
\raisebox{-1pt}{$\stackrel{\ast}{p}_i$}\, =\,
\gamma^{2i} p_i\, -\, (q-q^{-1})\sum_{k=1}^{i-1} \gamma^{2k} p_k
\raisebox{-1pt}{$\stackrel{\ast}{p}_{i-k}$}\, .
\ee
}\end{rem}

\begin{cor}
\lb{cor3.33}
Let $\R$ be a skew invertible {\rm R}-matrix of the $GL_q(n)$ type (in which case $\D$
is invertible, see proposition \ref{prop2.9}).
Then for the  matrix $M$ (\ref{M}) generating the RE algebra
${\goth L\scriptstyle GL_q(n)}[\Rd]\subset{\goth D\scriptstyle GL_q(n)}[\R,\gamma]$
following Cayley-Hamilton  identity is valid:
\be
\nn
%\lb{char-M}
\sum^{n}_{i=0}\, (-1/q)^{i}\, \stackrel{\ast}{a}_{i}\, M^{n-i}\, =\,
\sum^{n}_{i=0}\, ({-\gamma^2/ q})^i\, a_{i}\, M^{n-i}\, =\, 0 \, .
\ee
In the spectrally completed algebra
${\overline{\goth L}\scriptstyle GL_q(n)}[\Rd]
\subset{\overline{\goth D}\scriptstyle GL_q(n)}[\R,\gamma]$
this identity assumes a completely factorized form
\be
\lb{CH-factor-M}
\prod_{\alpha =1}^n\bigl(M - {\gamma^2\mu_{\alpha}\over q} I\bigr)\, =\, 0\, .
\ee
\end{cor}

With the factorized characteristic identity (\ref{CH-factor-M}) one can construct yet another resolution
of matrix unity (c.f. with (\ref{proj}))
\be
\lb{proj2}
S^\alpha\, :=\, \prod_{\beta=1\atop \beta\neq\alpha}^n{\bigl
(M - \gamma^2q^{-1}\mu_{\beta} I\bigr)\over \gamma^2q^{-1}(\mu_\alpha-\mu_\beta)}\; :\quad
S^\alpha S^\beta\, =\, \delta_{\alpha\beta}\, S^\alpha\, , \quad
\sum_{\alpha=1}^n S^\alpha\, =\, I\, ,
\ee
\vspace{-6mm}\ni
so that
\be
\nn
%\lb{M-P}
M\, S^\alpha\, =\, S^\alpha M\, =\, \gamma^2q^{-1}\mu_\alpha S^\alpha\, .
\ee
Relation between the two sets of projectors $P^\alpha$ and $S^\alpha$ is explained in
the following proposition.

\begin{prop}
\lb{prop3.30}
In the spectrally completed algebra
${\overline{\goth D}\scriptstyle GL_q(n)}[\R,\gamma]$
($\overline{{\goth D}}{\scriptstyle SL_q(n)}[\R]$) one has
\be
\lb{PTS}
P^\alpha T S^\beta \, =\, \delta_{\alpha\beta}\, P^\alpha T\qquad\mbox{or, equivalently,}\qquad
P^\alpha T\, =\, T S^\alpha
\, .
\ee
\end{prop}

\ni
{\em Proof.~}
Taking into account relations $TM=LT\,$, (\ref{L-P}) and (\ref{T-mu}) one
finds
$$
P^\alpha T M\, =\, P^\alpha L T\, =\, q\mu_\alpha (P^\alpha T)\, =\,
(P^\alpha T) \gamma^2 q^{-1} \mu_\alpha\, .
$$
Hence, in view of \eqref{proj2}
\ba
\nn
P^\alpha T S^\beta\, =\,
P^\alpha T\, \prod_{\sigma\neq \beta}
\frac{\gamma^2q^{-1}(\mu_\alpha -\mu_\sigma)}{\gamma^2q^{-1}
(\mu_\beta -\mu_\sigma)}\,.
\ea
In case  $\alpha\neq \beta$ the factor with $\sigma=\alpha$ in the product vanishes.
In case  $\alpha=\beta$ (and so, $\sigma\neq\alpha$) all the terms in the product are equal to 1.
So, the relation above reduces to the first equality in \eqref{PTS}.
\hfill$\blacksquare$

\subsection{Derivation of dynamical R-matrix}
\lb{subsec3.5}
In \cite{AF1} A.Alekseev and L.Faddeev used  dynamical R-matrix
in their construction of the Heisenberg double algebra.
Namely, they observed an appearance of the classical dynamical r-matrix
in the Poisson relations for certain classical variables
and then, by postulating a quantum counterpart of those relations,
they derived defining formulas (as in the definition \ref{3.3.1})
for the algebra ${\goth D\goth G}[\R,\gamma]$.

In this section we are aimed to
explain an origin of the dynamical R-matrix in the context of the HD algebras.
We show that the dynamical R-matrix -- $\R(\mu)_{\alpha\beta}$ --
appears in the permutation relations for matrix components of the matrices
\be
\lb{W}
W^\alpha\, :=\, P^\alpha T\, =\, T S^\alpha\, ,
\ee
and the arguments of the dynamical R-matrix are just the spectral variables
$\mu_{\alpha}$. In a sense, we solve an inverse problem to that considered in \cite{AF1}.

\smallskip
Recall the definition of two projectors associated with the Hecke type R-matrix
(see (\ref{H}))
\be
\lb{A-S}
A^{(2)}\,=\,
%\Rro(a^{(2)})\, =\,
{q\id - \R_1\over q+q^{-1}}\, ,\qquad
S^{(2)}\, =\, {q^{-1}\id + \R_1\over q+q^{-1}}\, .
\ee
These projectors, called the {\em antisymmetrizer} and the {\em symmetrizer},
serve for suitable separation of the different eigenspaces of $\R$.

\begin{theor}
\lb{theor3.35}
In the completed HD algebra
$\overline{{\goth D}}{\scriptstyle GL_q(n)}[\R,\gamma]$
($\overline{{\goth D}}{\scriptstyle SL_q(n)}[\R]$)
the matrices $W^\alpha$ (\ref{W}) satisfy  relations
\ba
\lb{AWWS}
S^{(2)}\Bigl\{ W^\alpha_1 W^\beta_2 +W^\beta_1 W^\alpha_2\Bigr\}A^{(2)}
= A^{(2)}\Bigl\{ W^\alpha_1 W^\beta_2 +W^\beta_1 W^\alpha_2\Bigr\} S^{(2)}
&=& 0\quad\forall\,\alpha,\beta ,
\\[2mm]
\lb{SWWS}
S^{(2)}\Bigl\{ (\mu_\beta - q^2\mu_\alpha)W^\alpha_1 W^\beta_2\, +\,
(\mu_\alpha - q^2\mu_\beta) W^\beta_1 W^\alpha_2\Bigr\}S^{(2)}
&=& 0\quad\forall\,\alpha\neq\beta ,
\\[2mm]
\nn
A^{(2)}\Bigl\{(\mu_\alpha - q^2\mu_\beta) W^\alpha_1 W^\beta_2\, +\,
(\mu_\beta - q^2\mu_\alpha)W^\beta_1 W^\alpha_2 \hspace{9.2mm}&&
\\
\lb{AWWA}
-\,(q^4-1)\,\mu_\alpha\,\varphi_{\alpha\beta}\,W^\alpha_1 W^\alpha_2\, -
\,(q^4-1)\,\mu_\beta\,\varphi_{\beta\alpha}\,W^\beta_1 W^\beta_2
\Bigr\}A^{(2)} &=& 0\quad\forall\,\alpha\neq\beta ,
\ea
where $\varphi_{\alpha\beta}:=\prod_{\sigma\neq\alpha ,\beta}
{\mu_\sigma - q^2\mu_\alpha\over \mu_\sigma -\mu_\beta}$. Relations (\ref{AWWS})--(\ref{SWWS})
and (\ref{T-mu}) (together with the appropriate conditions
on the spectral variables $\mu_\alpha$) define the algebra
$\overline{{\goth D}}{\scriptstyle GL_q(n)}[\R,\gamma]$
($\overline{{\goth D}}{\scriptstyle SL_q(n)}[\R]$)
in terms of generators $\,W^\alpha$, $\mu_\alpha$, $\;\alpha =1,\dots ,n$.
\end{theor}

\ni
{\em Proof.~}
Consider the product $W_1^\alpha W_2^\beta$, where $\alpha\neq\beta$. With the
help of (\ref{T-mu}) and (\ref{LkT}) we can reorder terms of the product in a following way
$$
W_1^\alpha W_2^\beta\, =\,
\frac{(L_1-q\mu_\beta)(L_{\overline{2}}J_2-q^3\mu_\alpha)}{q^2(\mu_\alpha -\mu_\beta)
(\mu_\beta -q^2\mu_\alpha)}\,W^{\alpha\beta}_{12},\quad\;\;
W^{\alpha\beta}_{12}:=\!\!
\prod_{\sigma\neq\alpha ,\beta}\!\!
\frac{(L_1-q\mu_\sigma)(L_{\overline{2}}J_2-q\mu_\sigma)}{q^2(\mu_\alpha -\mu_\sigma)
(\mu_\beta -\mu_\sigma)}\, T_1 T_2.
$$
Here factor $W^{\alpha\beta}_{12}$ commutes with the R-matrix
$\R_{12}$, which follows by the same arguments
 as in the proof of corollary \ref{cor3.25}, see below
(\ref{lhs_of}). We now extract symmetric and antisymmetric parts of the
product using projectors (\ref{A-S})
\ba
\nn
S^{(2)} W^\alpha_1 W^\beta_2&\!\!\!\! =\!\!\!\!& S^{(2)}\,
\frac{L_1 L_{\overline{2}}J_2 + q^4\mu_\alpha\mu_\beta\id -{q^2(\mu_\beta+\mu_\alpha)\over q+q^{-1}}
(L_1+L_{\overline{2}}J_2)+{\mu_\beta -q^2\mu_\alpha\over q+q^{-1}}
(L_1-L_{\overline{2}})}{q^2(\mu_\alpha -\mu_\beta)(\mu_\beta - q^2\mu_\alpha)}\,
W^{\alpha\beta}_{12} ,
\\[-2mm]
\lb{SWW}
&&
\\[1mm]
\nn
A^{(2)} W^\alpha_1 W^\beta_2&\!\!\!\! =\!\!\!\!& A^{(2)}\,
\frac{L_1 L_{\overline{2}}J_2 + q^4\mu_\alpha\mu_\beta\id -{\mu_\beta+q^4\mu_\alpha\over q+q^{-1}}
(L_1+L_{\overline{2}}J_2)+{q^2(\mu_\beta -q^2\mu_\alpha)\over q+q^{-1}}
(L_1-L_{\overline{2}})}{q^2(\mu_\alpha -\mu_\beta)(\mu_\beta - q^2\mu_\alpha)}\,
W^{\alpha\beta}_{12} .
\\[-2mm]
\lb{AWW}
\ea
Here we
separated linear in $L$ terms with the opposite symmetry properties
\be
\lb{L-AS}
\R_1 (L_1 + L_{\overline{2}}J_2) = (L_1 + L_{\overline{2}}J_2)R_1\, ,\qquad
\R_1 (L_1 - L_{\overline{2}}) = -(L_1 - L_{\overline{2}}) \R_1^{-1}\, ,
\ee
which was done  by the use of relation
$$
q^3 \mu_\alpha L_1 + q\mu_\beta L_{\overline{2}}J_2 =
\frac{q(\id +\R_1^{-2})}{(q+q^{-1})^2}\Bigl\{
(\mu_\beta \R_1^2 +q^2\mu_\alpha\id)(L_1+L_{\overline{2}}J_2) +
(q^2\mu_\alpha - \mu_\beta)(L_1 -L_{\overline{2}})
\Bigr\}\, .
$$
The symmetry properties (\ref{L-AS}) imply, in particular, that the only term
contributing to  expressions
$A^{(2)} W^\alpha_1 W^\beta_2 S^{(2)}$ and $S^{(2)} W^\alpha_1 W^\beta_2 A^{(2)}$
is the one proportional to $(L_1-L_{\overline{2}})$,
while the terms $(L_1 L_{\overline{2}}J_2)$, $I$ and $(L_1+L_{\overline{2}}J_2)$
contribute to
$S^{(2)} W^\alpha_1 W^\beta_2 S^{(2)}$ and $A^{(2)} W^\alpha_1 W^\beta_2 A^{(2)}$.

It is now straightforward to check that formulas (\ref{AWWS}) and (\ref{SWWS})
follow from relations (\ref{SWW}), (\ref{AWW}). To check formula (\ref{AWWA}) one needs
also similar expression for the product $W^\alpha_1 W^\alpha_2$:
$$
W^\alpha_1 W^\alpha_2\, =\,
\frac{L_1 L_{\overline{2}}J_2 + q^2\mu_\beta^2\id -q\mu_\beta
(L_1+L_{\overline{2}}J_2)}{q^2\varphi_{\alpha\beta}(\mu_\alpha -\mu_\beta)(q^2\mu_\alpha - \mu_\beta)}\,
W^{\alpha\beta}_{12}\, ,
$$
and analogous formula for $W^\beta_1 W^\beta_2$. Here factor $\varphi_{\alpha\beta}$
was defined in the proposition.

It lasts checking that defining relations for the algebras
$\overline{{\goth D}}{\scriptstyle GL_q(n)}[\R,\gamma]$
and $\overline{{\goth D}}{\scriptstyle SL_q(n)}[\R]$ can be derived from
(\ref{AWWS})--(\ref{AWWA}) and (\ref{T-mu}).
It is convenient to check relations for the matrices $T$ and $LT$:
\ba
\lb{T-LT}
\R_1 T_1 T_2 = T_1 T_2 \R_1, \;\;\;
\R_1 (LT)_1 (LT)_2 = (LT)_1 (LT)_2 \R_1,\;\;\;
\gamma^2  T_1 (LT)_2 = \R_1(LT)_1 T_2 \R_1 .\;\;\;
\ea
For  $GL_q(n)$ and $SL_q(n)$ types HD algebras, where $T$ is invertible, these formulas imply
(\ref{REA}) and (\ref{LT}). Substituting expressions
$$
T={\textstyle \sum_{\alpha =1}^n} W^\alpha,\qquad
LT = q\,{\textstyle \sum_{\alpha =1}^n}\,\mu_\alpha W^\alpha
$$
one  can easily prove that the first two relations (\ref{T-LT}) follow
from (\ref{AWWS}) and (\ref{T-mu}).
Checking the last formula in (\ref{T-LT}) is also straightforward, although more lengthy.
To this end one has to use the whole set of relations for $W$'s and to
take into account the identity $\sum_{\alpha\neq\beta}\varphi_{\alpha\beta}=1$ .
\hfill$\blacksquare$

\begin{cor}
\lb{cor3.36}
Relations (\ref{AWWS})--(\ref{AWWA}) can be equivalently written as
\ba
\lb{S-dyn}
S^{(2)}\Bigl[ W_1^\alpha W_2^\beta \R_{1}\, -\,
{\sum_{\alpha'\!,\beta' =1}^n} {\R^S(q;\mu)}_{\alpha' \beta'}^{\alpha\;\,\beta}\,
W_1^{\alpha'} W_2^{\beta'}\Bigr]&=&0\, ,
\\[2mm]
\lb{A-dyn}
A^{(2)}\Bigl[ W_1^\alpha W_2^\beta \R_{1}\, -\,
{\sum_{\alpha'\!,\beta' =1}^n} {\R^A(q;\mu)}_{\alpha' \beta'}^{\alpha\;\,\beta}\,
W_1^{\alpha'} W_2^{\beta'}\Bigr]&=&0\, ,
\ea
where $n^2\times n^2$ matrix $\R^{S}(q;\mu)$ has following nonzero components
$$
{\R^S\,}^{\alpha\alpha}_{\alpha\alpha}\, =\, q, \quad\;
{\R^S\,}^{\alpha\beta}_{\alpha\beta}\, =\, -{(q-q^{-1})\mu_\beta\over \mu_\alpha -\mu_\beta}, \quad\;
{\R^S\,}^{\alpha\beta}_{\beta\alpha}\, =\, {q^{-1}\mu_\alpha - q\mu_\beta\over \mu_\alpha -\mu_\beta}
\quad\forall\,\alpha\neq\beta,
$$
and $n^2\times n^2$ matrix $\R^{A}(q;\mu)$ has nonzero components at the same places
as $\R^{S}$ with values $\R^{A}(q;\mu)=\R^{S}(-q^{-1};\mu)$,
and the additional nonzero components
$$
{\R^A\,}^{\alpha\alpha}_{\beta\alpha}\, =\, -{\R^A\,}^{\alpha\alpha}_{\alpha\beta}\, =\,
{(q^4-1)\,\mu_\alpha\,\varphi_{\alpha\beta}\over q(\mu_\alpha -\mu_\beta)} \quad
\forall\,\alpha\neq\beta .
$$
Both matrices $\R^{S/A}(q;\mu)\equiv \R^{S/A}(\mu)$  satisfy
{\em dynamical Yang-Baxter equation}:
\be
\lb{DYBE}
\R\bigl(\mu\bigr)^{12}\R\bigl(\nabla^1(\mu)\bigr)^{23}\R\bigl(\mu\bigr)^{12}\, =\,
\R\bigl(\nabla^1(\mu)\bigr)^{23}\R\bigl(\mu\bigr)^{12}\R\bigl(\nabla^1(\mu)\bigr)^{23}\, .
\ee
Here superscript labels denote endomorphism spaces for the spectral indices, e.g.,
${\R(\mu)}^{\alpha_1\alpha_2}_{\beta_1\beta_2}\equiv \R(\mu)^{12}$, and $\nabla^1$
is a diagonal {\em finite shift operator}
\be
\lb{nabla-1}
\nabla^1\, =\, \mbox{\rm diag}\{ \nabla^\alpha \}_{\alpha =1}^{\phantom{\alpha =}n}\,:\quad
\nabla^\alpha(\mu_\beta)\,:=\, q^{2\delta_{\alpha\beta}}\gamma^{-2}\mu_\beta\, .
\ee
\end{cor}

\ni
{\em Proof.~}
Apply the symmetrizer $S^{(2)}$ and the antisymmetrizer $A^{(2)}$ from the right to both sides of
the equalities (\ref{S-dyn}), (\ref{A-dyn}). The resulting projections is easy to compare with
(\ref{AWWS})--(\ref{AWWA}).

To prove the dynamical Yang-Baxter equation
for the matrices $\R^A(q;\mu)$ and $\R^S(q;\mu)$
we consider, respectively, the following cubic terms:
$$
A^{(3)} W_1^\alpha W_2^\beta W_3^\sigma
\qquad \mbox{and}\qquad
S^{(3)} W_1^\alpha W_2^\beta W_3^\sigma\, .
$$
Here 3-antisymmetrizer $A^{(3)}=\Rro(a^{(3)})$  is the R-matrix realization of the idempotent $a^{(3)}$
(see (\ref{q-anti1}), (\ref{q-anti2})), and 3-symmetrizer $S^{(3)}$ is a similar projector
which differs from $A^{(3)}$ by substitution $q\leftrightarrow -q^{-1}$ in the formulas
(\ref{q-anti1}), (\ref{q-anti2}). Now applying
two equal operators $\R_1 \R_2 \R_1$ and $\R_2 \R_1 \R_2$ from the right side
to these terms and using relations (\ref{S-dyn}), (\ref{A-dyn}) and (\ref{T-mu})
one eventually proves (\ref{DYBE}) for $\R^{A/S}(q;\mu)$.
\hfill$\blacksquare$

\begin{rem}
{\rm\small
The dynamical R-matrix $\R^S(q;\mu)$ was constructed in \cite{F0,AF1,Is}.
A review on the dynamical Yang-Baxter  equation and the dynamical R-matrices
is given in \cite{ES}. It is  surprising that in our approach the dynamical R-matrices
$\R^{A/S}(q;\mu)$, being the solutions of the nonlinear finite difference equation (\ref{DYBE}),
are calculated by solving a system of (at most three) linear equations.
}
\end{rem}

In conclusion of the section we comment how relations (\ref{S-dyn}) can be reduced to
dynamical quadratic relations considered in \cite{F0,AF1}.
Recall that a (Hecke type) {\em quantum plane} ${\cal V}[\R]$
is an algebra generated by components of
vector $\{x_i\}_{i=1}^{\dim V}$ subject to relations
\be
\lb{V[R]}
x_{\langle 1|}x_{\langle 2|}\, A^{(2)}\, =\, 0\qquad\Leftrightarrow\qquad
x_{\langle 1|}x_{\langle 2|}\, S^{(2)}\, =\, x_{\langle 1|}x_{\langle 2|}\, .
\ee
In the tensor product algebra
${\cal V}[\R]\otimes \overline{{\goth D}}{\scriptstyle GL_q(n)}[\R,\gamma]$
consider a rectangular
%$\, n\times\dim V$
matrix
$$
\Lambda^\alpha_i\, :=\, {\textstyle \sum_{j=1}^{\dim V}}\,x_j\otimes W_{ji}^\alpha\, ,
\qquad\alpha =1,\dots n\, , \;\; i=1,\dots ,\dim V\, .
$$
As a consequence of (\ref{S-dyn}), (\ref{V[R]}) the matrix components of $\Lambda$
fulfill relations
\be
\lb{dyn-quadr}
\Lambda^{|1\rangle}_{\;\langle 1|}\,\Lambda^{|2\rangle}_{\;\langle 2|}\, \R_{12}\, =\,
\R^S(q;\mu)^{\!12}\, \Lambda^{|1\rangle}_{\;\langle 1|}\,\Lambda^{|2\rangle}_{\;\langle 2|}\, .
\ee
Assume additionally that {\em\, i)} $\dim V = n$,\; and {\em\, ii)} the quantum plane
admits a one dimensional representation $\chi:\; {\cal V}[\R]\rightarrow {\mathbb C}\,$
(note that both these conditions are satisfied for the R-matrices from the example
\ref{2.4.4}).

It is the square matrix
$\chi (\Lambda)\in\overline{{\goth D}}{\scriptstyle GL_q(n)}[\R,\gamma]$
whose dynamical quadratic relations (\ref{dyn-quadr}) were introduced in
\cite{F0,AF1} and also investigated  in \cite{HIOPT,FHIOPT}.

\section{Discrete time evolution on quantum group cotangent bundle}
\lb{sec4}

\subsection{Automorphisms of the Heisenberg double algebra}
\lb{subsec4.1}

In this section we investigate a
sequence of automorphisms on the HD algebra $ {\goth D\goth G}[\R,\gamma]$.
These automorphisms were introduced by A.\,Alekseev and L.\,Faddeev \cite{AF1,AF2}, who
interpreted them as a discrete time evolution of a  $q$-deformed quantum isotropic Euler top.
The automorphisms $\theta^k: {\goth D\goth G}[\R,\gamma]\rightarrow{\goth D\goth G}[\R,\gamma]$
are given on generators
\ba
\nn
&\{ T, \; L \} \stackrel{\theta^{k}}{\longrightarrow} \{ T(k) , \; L(k) \} \, ,
\quad \forall\, k=0,1,2,\dots\, ,&
\\
\lb{aut1}
&T(0) := T \; , \quad T(k+1): = L T(k) = L^{k+1} \, T \; , \quad L(k) := L \, .&
\ea
It is easy to see (c.f. \eqref{T-LT}) that the map $\theta$ agrees with the defining relations (\ref{RTT}), (\ref{REA}),
(\ref{LT}) of the algebra ${\goth D\goth G}[\R,\gamma]$.
Less obvious is its consistency with the $SL_q(n)$ type reduction conditions.
\begin{prop}
\lb{prop4.0.1}
The map $\theta$ (\ref{aut1}) defines an automorphism of the algebra
${\goth D\scriptstyle SL_q(n)}[\R]$.
\end{prop}

\ni
{\em Proof.~} It is necessary to check that  $\Rdet (LT) =1$ in the $SL_q(n)$ case.
To this end, we use formula
\ba
\nn
(LT)_1(LT)_2\dots (LT)_k &=&
\gamma^{-k(k-1)} Z_k \bigl( L_{\overline{1}}L_{\overline{2}}\dots L_{\overline{k}}\bigr)
\bigl(T_1 T_2\dots T_k \bigr)
\ea
to separate matrices $L$ and $T$ in the expression for $\Rdet (LT)$. This formula follows
from (\ref{LT}), (\ref{REA-forms}), (\ref{LR-comm}) and (\ref{J-k}) by induction on $k$.

The calculation of $\Rdet (LT)$ proceeds as follows
\ba
\nonumber
{\det}_R LT &\!\!\!\!\!:=&\!\!\!\! \Tr{1,\dots, n}\left(A^{(n)} (LT)_1\dots(LT)_n\right) =
\gamma^{n(n-1)}\Tr{1,\dots, n}\left(A^{(n)} Z^{(n)}L_{\overline{1}}\dots L_{\overline{n}}
T_1\dots T_n\right)
\\[1mm]
\nn
&\!\!\!\!\!=&\!\!\!\!\!(\gamma q)^{-n(n-1)}\,\Tr{1,\dots, n}\bigl(A^{(n)}L_{\overline{1}}\dots L_{\overline{n}}\bigr)
\,\Tr{1,\dots, n}\bigl(A^{(n)}T_1\dots T_n\bigr)
\\[1mm]
\lb{evol-detT}
&\!\!\!\!\!=&\!\!\!\!\!q^n\gamma^{-n(n-1)}\,\RTr{1,\dots, n}\bigl(A^{(n)}L_{\overline{1}}\dots L_{\overline{n}}\bigr)\,
\Rdet T\, =\, \bigl(q\gamma^{-n}\bigr)^{n-1}\, q\, a_n\, \Rdet T\, ,
\ea
and so, under conditions $\Rdet T=1$,~ $a_n =q^{-1} 1$,~ $\gamma^n=q$~
we have $\Rdet (LT)=1$.
Here in the second line we substituted $A^{(n)}Z_n=q^{-n(n-1)}A^{(n)}$ and used the
condition  ${\rm rk}A^{(n)} = 1$; in the last line we applied (\ref{D-A}) and the definitions of
$\Rdet T$ and $a_n$.
\hfill$\blacksquare$

\smallskip
In what follows we will investigate the automorphisms (\ref{aut1}) for
HD algebras of the types ${\goth D\scriptstyle GL_q(n)}[\R,\gamma]$ and
${\goth D\scriptstyle SL_q(n)}[\R]$.
A key point for their dynamical interpretation
is a possibility to write down following ansatz
\be
\lb{aut2}
T(k+1)\, =\, L\, T(k)\, =\, \bigl(q a_n\bigr)^{1/n}\,\Theta \, T(k) \, {\Theta}^{-1} \, ,\quad
\mbox{where}\;\; \Theta\in{\overline{\goth C\goth h}}[\R]\, .
\ee
Here the dynamical process -- evolution -- is thought as an inner HD
algebra automorphism, and
$\Theta$ plays a role of the {\em evolution operator}.
As the evolution keeps $L$ unchanged,
it is natural to assume that $\Theta$ belongs to the center
of the RE algebra generated by the
matrix $L$. More specifically, we will look for $\Theta$
as a formal power series in spectral variables $\mu_\alpha$, $\alpha=1,\dots ,n$,
which we denote as
${\overline{\goth C\goth h}}[\R]$.
We also note that the condition $\Theta\in{\overline{\goth C\goth h}}[\R]$ makes the ansatz
manifestly covariant with respect to both left and right coactions (\ref{delta-l}), (\ref{delta-r}).

Factor $(q a_n)^{1/n}$ in the ansatz (\ref{aut2}) becomes trivial for the $SL_q(n)$ type HD algebra.
In the $GL_q(n)$ case one adds this scaling factor to make the ansatz
consistent with the evolution of $\Rdet T$, see (\ref{evol-detT}). One assumes
following relation for the newly introduced element $a_n^{1/n}$ (c.f. with (\ref{T-an}))
\be
\lb{T-a_n}
T\, a_n^{1/n}\, =\, \bigl(q\gamma^{-n}\bigr)^{2/n}\, a_n^{1/n}\, T\, .
\ee
Then, consistency of (\ref{aut2}) and (\ref{evol-detT}) results in a commutativity of
$\Rdet T$ with $\Theta$:
\be
\lb{theta-detT}
\Rdet T\, \Theta\, =\, \Theta\, \Rdet T\, ,
\ee
which again trivializes in the $SL_q(n)$ case.

\begin{rem}{\rm\small
The action of the automorphisms $\theta^k$ on $T$
can be equally treated as write multiplications by powers of the left invariant
matrix $M$:
$$
T(k+1)\, =\, T(k)\, M\, =\, T\, M^k\, , \quad M(k)\, =\, M\, .
$$
The relation (\ref{M}) between $L$ and $M$ would no more be valid if one would treat them as
{\em quantized right and left invariant Lie derivatives} acting on {\em quantized external
algebra of differential forms over matrix group}. In this case one would have
a two-parametric series of automorphisms:
$$
\{ T, \; L,\; M \}\; \stackrel{\theta^{(k,m)}}{\longrightarrow}\; \{ L^kTM^m , \; L, \; M \} \, ,
\quad \forall\, k\geq 0,\; m\geq 0\, .
$$
}\end{rem}

\begin{example}
\lb{example-ribbon}
{\rm %\small
Let us show that in the ribbon Hopf algebra setting
the ribbon element $\upsilon\in {\goth A}_{\cal R}$ generates the evolution (\ref{aut2})
in the smash product algebra ${\goth A}_{\cal R}\sharp\, {\goth A}_{\cal R}^*$.
For this we first have to specify pairing for the ribbon element.
Using the definition (\ref{ribbon-def})  and relations (\ref{CD-iz-univ}), (\ref{C-D}) and setting
$\eta =q^{1/n}$ as in the example \ref{2.4.4} we calculate
\be
\nn
\langle T,\, \upsilon^2\rangle\, =\,
\langle T,\, u S(u)\rangle\, =\, \rho_V(u)\, \rho_V(S(u))\, =\,
\eta^2 \D \,\C\, =\, q^{2({1\over n}-n)}\,\id\, .
\ee
Therefore, taking into account centrality of the ribbon element $\upsilon$ in  ${\goth A}_{\cal R}$,
it is  natural to define
\be
\nn
\langle T,\, \upsilon\rangle\, =\, q^{({1\over n}-n)}\,\id\, .
\ee
Using this formula and relations (\ref{ribbon-def}), (\ref{L-iz-univ}), \eqref{smash}
we now calculate conjugation of matrix $T$ with the ribbon element
\ba
\nn
\upsilon\, T \, \upsilon^{-1} &=&
(\upsilon\otimes id)\,
\langle  id \otimes T,\, \Delta(\upsilon^{-1}) \rangle\, T\, =\,
(\upsilon\otimes id)\,
\langle id \otimes T,\,(\upsilon^{-1}\otimes \upsilon^{-1})
{\cal R}_{21}{\cal R}_{12}\rangle\, T
\\[1mm]
\nn
&=&
\langle T,\, \upsilon^{-1}\rangle\,
\langle id \otimes T , \,
{\cal R}_{21} {\cal R}_{12} \rangle\, T\, =\, LT\, .
\ea
Note that defining relations for the evolution operator $\Theta$
(as a function of the spectral variables
$\mu_\alpha$)
and for the ribbon element $\upsilon$, both admit multiple solutions.\,\footnote{
The ribbon element  is defined modulo central factor $z\in {\goth A}_{\cal R}$:
$z^2 =1 \, , \;\; S(z) = z \, , \;\; \epsilon(z) =1 \, , \;\; \Delta(z) = z \otimes z \,$.
For the evolution operator $\Theta(\mu)$ different solutions are constructed in the next sections.
}
Therefore, a problem of finding explicit expression of the ribbon element $\upsilon$
in terms of spectral variables $\mu_\alpha$ demands further investigations.
}
\end{example}

\subsection{Equations for the evolution operator $\Theta$}
\lb{subsec4.2}

Using the results of section \ref{sec3} it is straightforward to derive
equations for $\Theta$.
We consider in details, the evolution in the $SL_q(n)$ type HD algebra.
In this case we assume
\be
\lb{cond-sl}
\Theta=\Theta(\mu_1,\dots ,\mu_n)\, , \qquad\mbox{where}\quad
a_n={\textstyle \prod_{\alpha =1}^n}\,\mu_\alpha =q^{-1} \quad\mbox{and}\quad
\gamma=q^{-1/n}\, .
\ee
Applying from the left the projector $P^\alpha$ to  both sides of
(\ref{aut2}) we obtain
$$
q\mu_\alpha \, (P^\alpha T)\, =\, \Theta (P^\alpha T) \Theta^{-1}\, ,
\quad \forall\, \alpha =1,\dots ,n\, .
$$
Multiplying this equality by $\Theta$ from the right and permuting $\Theta$ with $P^\alpha T$
in the left hand side with the help of (\ref{T-mu}) we finally get
$$
q\mu_\alpha\, \Theta(q^{-2/n}\mu_1, \dots ,q^{2-2/n}\mu_\alpha, \dots ,q^{-2/n}\mu_n)\,
(P^\alpha T)\, =\, \Theta(\mu_1,\dots ,\mu_n)\, (P^\alpha T)\, .
$$
We state the result in a following proposition

\begin{prop}
\lb{prop4.3}
For the Heisenberg double algebra ${\goth D\scriptstyle SL_q(n)}[\R]$ the evolution
operator $\Theta(\mu_\alpha)$ in (\ref{aut2}), (\ref{cond-sl}) is a solution of equations
\be
\lb{Sl-evol}
q\mu_\alpha\, \Theta\bigl(\nabla^\alpha(\mu_\beta)\bigr)\,
=\, \Theta(\mu_\beta)\quad\forall\,\alpha =1,\dots ,n\, ,
\ee
where $\nabla^\alpha$ are finite shift operators introduced in (\ref{nabla-1}). In the
$SL_q(n)$ case their actions are
\be
\lb{nabla}
\nabla^\alpha(\mu_\beta)\, :=\, q^{2 X_{\alpha\beta}}\,\mu_\beta\, ,
\qquad
X_{\alpha\beta}\, :=\, \delta_{\alpha\beta} \, -\, {1\over n}\, \quad \forall\,
\alpha, \beta=1,\dots ,n\, .
\ee
The $n\times n$ matrix $X$ is a Gram matrix for the set of vectors
$\vec{e}^{\,\,*}_\alpha\in {\mathbb Q}^n$,
$\,\alpha = 1,\dots ,n\,$:
$$
\vec{e}^{\,\,*}_\alpha := {1\over n}\,
(\,\underbrace{-1, \dots, -1}_{\scriptstyle (\alpha-1)\;\mbox{\small times}},n-\! 1,-1,\dots ,-1\,)\, ,
\qquad
X_{\alpha \beta}\, =\, \langle\vec{e}^{\,\,*}_\alpha , \vec{e}^{\,\,*}_\beta\rangle\, .\vspace{-1mm}
$$
$X$ is positive semi-definite of the rank $\;n-1$ ($\,\,\sum_{\alpha =1}^n\vec{e}^{\,\,*}_\alpha =0\,$).

\smallskip
For the Heisenberg double algebra ${\goth D\scriptstyle GL_q(n)}[\R,\gamma]$ the evolution
operator $\Theta$ is suitably parameterized by variables $z:=(qa_n)^{1/n}$ and $\nu_\alpha$
\be
\nn
%\lb{nus}
\nu_\alpha := \mu_\alpha (q a_n)^{-1/n}\, ,\qquad \mbox{such that}\qquad
{\textstyle \prod_{\alpha =1}^n}\, \nu_\alpha\, =\, q^{-1}\,, \quad
\nu_\alpha \, \Rdet T\, =\, \Rdet T\, \nu_\alpha\;\; \forall\, \alpha\, .
\ee
The evolution equations for $\Theta(\nu_1,\dots ,\nu_n; z)$ read
\be
\lb{Gl-evol}
q\nu_\alpha\, \Theta(\nabla^\alpha(\nu_\beta);\, (q\gamma^{-n})^{2/n}\, z)\, =\,
\Theta(\nu_\alpha; z)\quad\forall\,\alpha =1,\dots ,n\, ,
\ee
where shift operators $\nabla^\alpha$ are defined as in (\ref{nabla}).
Since $\prod_{\alpha =1}^n \nabla^\alpha =1$, this system
is consistent provided that
$\,\Theta(\nu_\beta; \, q^2\gamma^{-2n}\, z)\,=\, \Theta(\nu_\beta; z)\,$
(c.f. with (\ref{theta-detT})). Demanding that $\Theta$ does not actually depend on $z$
one reduces (\ref{Gl-evol}) to (\ref{Sl-evol}).
\end{prop}

\ni
{\em Proof.~}
The $SL_q(n)$ case is already considered.
Taking into account relations
(\ref{T-a_n}) and (\ref{detT-mu})
a derivation of the evolution equations
in the $GL_q(n)$ case is the same.
\hfill$\blacksquare$

\medskip
In the next two subsections we will construct particular solutions of the
$SL$-type evolution equations (\ref{Sl-evol}).

\subsection{Solution in case $|q|<1$}
\lb{subsec4.3}

Let us look for solution of (\ref{Sl-evol}) as a series in
$\mu_\alpha$. Taking into account condition (\ref{cond-sl})
we exclude one dependent
variable, say $\mu_n\,$, from the expansion
\be
\lb{i.8}
\Theta(\mu_\alpha) =
\sum_{\;\; \vec{k} \in {\mathbb Z}^{^{n-1}}} \,
c(\vec{k}) \,
\mu_{1}^{k_1} \,  \mu_{2}^{k_2}  \dots  \mu_{n-1}^{k_{n-1}} \; .
\ee
where
${\mathbb Z}^{n-1} = \{(k_1, \dots,k_{n-1}): \; k_i \in {\mathbb Z} \}$,
and the coefficients $c(\vec{k})$ are  $\mathbb C$-valued functions on ${\mathbb Z}^{n-1}$.
Substitution of (\ref{i.8}) into  (\ref{Sl-evol}) gives conditions on the coefficients:
\be
\lb{i.8'}
c(\vec{k} + \vec{\epsilon}_\alpha)\, =\,
q^{\bigl(1+{2\over n} \sum_{\beta =1}^{n-1} A^*_{\alpha\beta} k_\beta\bigr) } \, c(\vec{k})
\quad \forall\,\alpha=1,\dots ,n-1.
\ee
Here
$
\vec{\epsilon}_\alpha := (\,\underbrace{0, \dots, 0}_{\lefteqn{\scriptstyle (\alpha-1)\;\mbox{\small times}}},1,0, \dots , 0)$,
and
$A^*_{\alpha\beta} := n\, X_{\alpha\beta}
%=(n\delta_{\alpha\beta} -1)\,
$ is a $(n-1)\times(n-1)$\vspace{-4mm} positive-definite
matrix.

\medskip
The general solution to (\ref{i.8'}) is
\be
\lb{i.12}
c(\vec{k})\, =\,
q^{{1\over n}\bigl\{ (\vec{k} ,\, A^* \, \vec{k})
+(\vec{1}\, \vec{k})\bigr\} } \,,\vspace{2mm}
\ee
where we choose normalization $c(\vec{0})=1$ and use notation
$\;(\vec{k} ,\, A^* \, \vec{k})=\sum_{\alpha,\beta=1}^{n-1}k_{\alpha} A^*_{\alpha\beta} k_\beta,\;$
and  $\;\vec{1}=(1,\dots ,1)\,$, so that
$\;(\vec{1},\, \vec{k})=\sum_{\alpha=1}^{n-1}k_\alpha\,$.

\begin{rem}{\rm\small
The matrix $A^*_{\alpha\beta}$   is a
Gram matrix of a lattice $A^*_{n-1}$
dual to the root lattice $A_{n-1}$ (see \cite{CoSl}, chapter 4, section 6.6).
The corresponding  quadratic form $(\vec{k},\, A^*\, \vec{k})$
is often referred to as Voronoi's principal form of the first type.
}
\end{rem}

The ansatz (\ref{i.8}) gives a particular solution
(\ref{i.12}) of the evolution equations, we  denote it $\Theta^{(1)}$.
Introducing a parameterization
\be
\lb{add}
q= \exp ( 2 \pi {\rm i\, } \tau)  , \quad
q^{1/n}\mu_\alpha =  \exp ( 2 \pi {\rm i \,} z_\alpha )  :\;
\sum_{\alpha =1}^{n-1} z_\alpha=0
,\quad
%(\vec{k} , \, \vec{z}) = \sum_{\alpha=1}^{n-1} k_\alpha z_\alpha \, ,\quad
\Omega_{\alpha\beta}={2\tau\over n}A^*_{\alpha\beta}\, =\, 2\tau\, (\delta_{\alpha\beta}-{1\over n})\, ,
\ee
we can write $\Theta^{(1)}$ as a Riemann theta function $\theta(\vec{z},\Omega)$ (see \cite{Mum,Ig})
\be
\lb{i.14}
\Theta^{(1)} (\mu_\alpha)\, =\, \theta(\vec{z}, \Omega)\, =\!\!\!
 \sum_{\;\; \vec{k} \in {\mathbb Z}^{^{n-1}}} \!\!
\exp \left\{ \pi  {\rm i \,} (\vec{k}, \, \Omega \, \vec{k})
+ 2 \pi {\rm i \,} (\vec{k} \, , \, \vec{z})
\right\}  .
\ee
Here $\tau$ is a {\em modular parameter}  and $\Omega$ is a {\em matrix of periods}.
Expression (\ref{i.14}) converges  either if $|q| < 1$, or if $q$ is a rational root
of unity, in which case the series can be truncated.

\begin{rem}
\lb{rem4.5}
{\rm\small
One can  present  formula (\ref{i.14}) in a manifestly covariant form:
\be
\nn
%\lb{i.14a}
\Theta^{(1)}\, \equiv\,
\Theta^{(1)}(\vec{z}, A^{*}_{n-1}, \tau)\, =
\!\!\! \sum_{\;\; \vec{\rm \, k} \in A^*_{n-1}} \!\!
\exp \left\{ 2 \pi {\rm i \,} \frac{\tau}{n} \, \langle\vec{\rm\, k},  \vec{\rm k}\rangle
+ 2 \pi {\rm i \,} \langle\vec{\rm k} \, ,  \vec{\rm \, z}\rangle
\right\} .\vspace{-3mm}
\ee
Here vectors $\vec{\rm\, k}=\sum_{\alpha=1}^{n-1} k_\alpha e^*_\alpha$
label vertices of the lattice $A^{*}_{n-1}$, and
$\;\vec{\rm \, z}=\sum_{\alpha=1}^{n-1} z_\alpha e_\alpha$, where
$e_\alpha= \epsilon_\alpha -\epsilon_n$, $\,\alpha=1,\dots ,n-1$, (see line below (\ref{i.8'}))
are basic vectors of the root lattice $A_{n-1}$:
$\langle e^*_\alpha,\, e_\beta\rangle=\delta_{\alpha\beta}$.
}
\end{rem}

In the simplest $SL_q(2)$ case the evolution operator $\Theta^{(1)}$ becomes the
Jacobi theta function:
\ba
\nn
\Theta^{(1)} (\mu_1)\, =\,
\sum_{k\in {\mathbb Z}} q^{\frac{1}{2} k(k+1)}\mu_{1}^{k}\, =\,
\sum_{k\in {\mathbb Z}} \exp(\pi {\rm i\,} k^2 \tau + 2\pi {\rm i\,} k z_1)\, =\,
\theta_3(z_1;\, q)\, ,
\ea
or, in a multiplicative form
\ba
\nn
\Theta^{(1)} (\mu_1)\, =\,
\prod_{n=1}^{\infty} (1-q^n)( 1 + q^n \mu_1 )(1+q^{n-1}/\mu_1) \, .
\ea

\subsection{Solution for arbitrary $q$ }
\lb{subsec4.4}

In this section we derive yet another particular solution of the evolution equations (\ref{Sl-evol}),
the one which is well defined for arbitrary values of $q$.
The idea of such solution was
proposed in  L.D. Faddeev's lectures on two-dimensional integrable quantum field theory
\cite{Fad} (see also \cite{Faddd}).
We use heuristic arguments inspired by
considerations in \cite{AF1}. For the moment we assume $\dim V =n$, so that the range of the indices
$\alpha$ and $i, j$ in the projectors $P^\alpha_{ij}$, $S^\alpha_{ij}$ is the same.
Consider following $n\times n$ matrices
\be
\nn
%\lb{UV}
U_{ij}\, :=\, \sum_{k=1}^n u_{ik}\, P^{\alpha =i}_{kj}\, ,\qquad
V_{ij}\, :=\, \sum_{k=1}^n  S^{\alpha =j}_{ik}\, v_{kj}\, ,
\ee
where the only restriction for the auxiliary parameters $u_{ij}$ and $v_{ij}$ is their
commutativity with the spectral variables $\mu_\alpha$
$$
[u_{ij}, \, \mu_\alpha]\, =\, [v_{ij}, \, \mu_\alpha]\, =\, 0\quad
\forall\, i,j,\alpha\, .
$$
As a result of the Cayley-Hamilton identities (\ref{CH-factor}), (\ref{CH-factor-M})
we have matrix relations
$$
UL\, =\, q D U\, , \qquad M V\, =\, {\gamma^2q^{-1}}\, V D\, ,\qquad \mbox{where}\;\;\;
D:=\mbox{diag}\{\mu_1, \dots ,\mu_n\}\, .
$$
Moreover, by (\ref{PTS}), matrix $Q:= UTV$ is diagonal
$$
Q = \mbox{\rm diag}\{w_1,\dots ,w_n\}\, , \qquad\mbox{where}\;\;
w_i\, :=\, \Bigl(u P^i T v\Bigr)_{ii}\, ,
$$
and, by (\ref{T-mu}), $w_i$ satisfy following permutation relations with $\mu_j$ and $z_j$:
\be
\lb{w-mu}
w_i\, \mu_j \, =\, q^{2\delta_{ij}} \gamma^{-2}\, \mu_j\, w_i\quad\Leftrightarrow\quad
w_i\, z_j\, =\, \Bigl(z_j + 2\tau (\delta_{ij}-1/n)\Bigr) w_i\, ,
\ee
where in the latter formula we used the $SL_q(n)$ type condition $\gamma=q^{1/n}$.

Assuming invertibility of the matrices $U$ and $V$
we can write diagonal decompositions for the matrices $L$, $M$ and $T$
$$
L\, =\, q\, U D U^{-1}\, ,\qquad
M\, =\,\gamma^2 q^{-1}\, V D V^{-1}\, ,\qquad
T\, =\, U^{-1} Q V^{-1}\, ,
$$
which after substitution into the ansatz (\ref{aut2}) reduce the evolution
equations to a following form
\be
\nn
%\lb{evol-new}
q D Q\, =\, \Theta\, Q\, \Theta^{-1}\qquad\Leftrightarrow\qquad
q \mu_i w_i \, =\, \Theta\, w_i\,\Theta^{-1}\, .
\ee
Taking into account (\ref{w-mu}) these equations clearly have following solution
\be
\lb{i.x-0}
\Theta^{(2)}(z_\alpha)\, :=\,
\exp \Bigl( -\frac{\pi {\rm i}}{2 \tau} {\textstyle \sum_{\beta=1}^{n} z_{\beta}^{2}} \Bigr)\, .
\ee
Now, it is easy to check that the function $\Theta^{(2)}$ fulfills the evolution equations
(\ref{Sl-evol}) without additional assumptions we made for the derivation.
Written in the independent variables
$\vec{z}=\{z_1,\dots ,z_{p-1}\}$ it reads
\ba
\lb{i.x}
\Theta^{(2)}(\vec{z})
\, =\, \exp \Bigl( -\frac{\pi {\rm i}}{\tau}
\sum_{1\leq \alpha\leq\beta\leq n-1} \!\!\!\!\!\!z_\alpha z_\beta \Bigr)
\, =\,
\exp \Bigl\{ -  \pi {\rm i\, } ( \vec{z},\, \Omega^{-1} \vec{z}) \Bigr\}\, ,
\ea
where the inverse matrix of periods is
\be
\nn
%\lb{gram3}
\Omega^{-1}_{\alpha\beta} = \frac{1}{2 \tau} \, (\delta_{\alpha\beta} + 1)
 = \frac{1}{2 \tau} \, A_{\alpha\beta} \; ,
\ee
and $A_{\alpha\beta} = \langle e_\alpha , e_\beta\rangle$
is the Gram matrix for the  root lattice $A_{n-1}$
(see remark \ref{rem4.5}).
Let us stress that
the logarithmic change of variables $\mu_\alpha \mapsto z_\alpha$ (\ref{add})
which was rather superficial in case of $\Theta^{(1)}$,
is inevitable for the derivation of $\Theta^{(2)}$.
\smallskip

Finally, we comment on relation between the two evolution operators
$\Theta^{(1)}=\theta(\vec{z},\,\Omega)$ and $\Theta^{(2)}$. The relation is based on a
functional equation for Riemann theta function
\ba
\nn
%\lb{multi-Jacobi}
\theta(\Omega^{-1}\vec{z},\, -\Omega^{-1})\, =\,
\Bigl(\det\bigl(\Omega/{\rm i}\bigr)\Bigr)^{1\over 2}\exp\Bigl\{
\pi {\rm i}(\vec{z},\, \Omega^{-1}\vec{z})\Bigr\}\,\theta(\vec{z},\, \Omega)\, ,
\ea
which is the special case of a more general modular functional equation
% in turn follows from the Poisson's summation formula
(for derivation and generalization see \cite{Mum}, chapter 2, section 5).
With our particular  matrix of periods $\Omega$ (\ref{i.14})
we  find
\be
\lb{Theta-1-2}
\Theta^{(2)}(\vec{z})\, =\, {1\over\sqrt{n}}
\Bigl({2\tau\over {\rm i}}\Bigr)^{n-1\over 2}
\frac{\theta(\vec{z},\, \Omega)}{\theta(\Omega^{-1}\vec{z},\, -\Omega^{-1})}\, .
\ee
Note that theta function in the denominator --
$\;\theta(\Omega^{-1}\vec{z},\, -\Omega^{-1})$ --
commutes with the elements of
${\goth D\scriptstyle SL_q(n)}[\R]$
and can be thought as an evolution operator on a `modular dual' quantum cotangent bundle
\cite{{Fadub}}.

% Then, we define the modular dual theta-function
% $$
% \Theta(\vec{z}, -\frac{1}{\Omega}) =
% \sum_{\vec{k} \in {\bf Z}^{n-1}} \exp \left\{ - \pi {\rm i}
% (\vec{k} \; \Omega^{-1} \vec{k} ) +2 \pi {\rm i} (\vec{k}, \; \vec{z}) \right\} =
% $$
% $$
% =  \sum_{\;\; {\bf l} \in \Lambda} \!\!
% \exp \left\{ - \frac{\pi {\rm i}}{2  \tau} \, ({\bf l}, \, {\bf l})
% + 2 \pi i \, ({\bf l} \, , \, \tilde{z})
% \right\} \; ,
% $$
% which serves a discrete evolution on a modular dual quantum cotangent bundle.

\appendix

\numberwithin{equation}{section}

\section{Pairing between spectral variables and  quantized functions}
\label{append-A}
Here we calculate pairing of the elementary symmetric functions $a_i$ (\ref{sigma}) with the
generators  of quantized functions $T^i_j$. We assume that $T$ and $a_i$
are realized respectively,  as elements of dual quasi-triangular Hopf algebras
${\goth A}^*_{\cal R}$ and ${\goth A}_{\cal R}$. We further extend
this pairing also for the spectral variables $\mu_\alpha$.

For the calculation we use formula
\be
\lb{pair-TL}
\langle T_1, \, L_2\rangle\, =\,\eta^{-2} q^{(n-{1\over n})}\R_{12}^2
\ee
which follows from the definitions (\ref{R-iz-univ}), (\ref{motiv1}), (\ref{dop1}),
(\ref{L-iz-univ}).

\begin{prop}
\lb{propA}
Let $a_i$ (\ref{dop1}), (\ref{L-iz-univ}) and $T$ (\ref{motiv1}) be elements of the
dual quasi-triangular Hopf algebras, respectively,
${\goth A}_{\cal R}$ and ${\goth A}^*_{\cal R}$
Assume that
the R-matrix $\R$ (\ref{R-iz-univ}) is $GL_q(n)$ type with scaling parameter
$\eta =q^{1/n}$ as in the example \ref{2.4.4}. Then
\be
\lb{pair-Ta}
\langle T, \, a_i\rangle\, =\, q^{-{3i/ n}}\,n_q^{-1} {n\choose i}_{\!\!q}
\Bigl\{n_q\, +\, q^{n+1}\, -\, q^{n-2i+1}\Bigr\}\id .
\ee
\end{prop}

\ni
{\em Proof.~} The calculation proceeds as follows
\ba
\nn
\langle T_1, \, a_i\rangle &=&
\langle T_1,\, \RTr{2,\dots ,i+1}\bigl( A^{(i)}
L_1\dots L_{\overline{i}}\bigr)^{\uparrow 1}\rangle
\\[1mm]
\nn
&=&q^{i(n-{3\over n})}\,
\RTr{2,\dots ,i+1}\bigl( A^{(i)\uparrow 1}(J_2\dots J_{i+1})
(J_1^{-1}\dots J_i^{-1})^{\uparrow 1}\bigr)
\\[1mm]
\nn
&=&q^{i(n-{3\over n})}\,
\RTr{2,\dots ,i+1}\bigl(\R_1\dots \R_i\, A^{(i)}\R_i\dots \R_1\bigr)
\\[1mm]
\nn
&=&
q^{i(n-{3\over n})}\,\Bigl\{
q^{-n}\,(n+1-i)_q\, i_q^{-1}\,\,
\RTr{2,\dots ,i}\bigl(\R_1\dots \R_{i-1}\, A^{(i-1)}\R_{i-1}\dots \R_1\bigr)
\\
\nn
&&
+\,(q-q^{-1})\,q^{-(n+2)(i-1)}\, n_q^{-1}\,{\textstyle {n\choose i}_q}\,\id_1
\Bigr\}
\\[1mm]
\nn
\dots &=&
q^{i(n-{3\over n})}\,{\textstyle {n\choose i}_q}\,
\Bigl\{ q^{-in}+
(q-q^{-1})\,q^{-(n+2)(i-1)}\, n_q^{-1}\,(1+q^2+\dots +q^{2(i-1)})\Bigr\}\id_1 .
%\\[1mm]
%\nn
%&=&
%q^{i(n-{3\over n})}\,n_q^{-1}\,{\textstyle {n\choose i}_q}
%\Bigl\{q^{-in}n_q\, +\, (q-q^{-1})\,q^{-(n+1)(i-1)}\,i_q\Bigr\}\id .
\ea
Here in the first line we substituted expression  for $a_i$ similar to
(\ref{a-altern}). In the second line we evaluated pairing using formulas
$\langle T_1,\, (L_{\overline{i}})^{\uparrow 1}\rangle = \eta^{-2} q^{(n-{1\over n})}J_{i+1}(J_i^{-1})^{\uparrow 1}$
following from (\ref{pair-TL}). In the third line we first, used the cyclic property
of the R-trace to evaluate the term
$(J_1^{-1}\dots J_i^{-1})^{\uparrow 1}$ on $(A^{(i)})^{\uparrow 1}$ and then, rearranged
the product $\,J_2\dots J_{i+1} = Z_{i+1} = J^\dagger_2\dots J^\dagger_{i+1}$,\, where
$\,J^\dagger_1:=\id, \;\; J^\dagger_{k+1}:=\R_{i-k+1}J^\dagger_k\R_{i-k+1}$,\, and evaluated
the term $J^\dagger_2\dots J^\dagger_{i}$ on $(A^{(i)})^{\uparrow 1}$. After that we
recollected terms in the product:
$(A^{(i)})^{\uparrow 1}J^\dagger_{i+1} = \R_1\dots \R_i\, A^{(i)}\R_i\dots \R_1$.
In the forth line we substituted $\R_i = \R_i^{-1} + (q-q^{-1})\id$ for one of $\R_i$s
and used formulas (\ref{dmat1}) and (\ref{dmat4}) to evaluate $\RTr{i+1}$. Then,
in the summand which is proportional to $(q-q^{-1})$ all the R-traces can be evaluated
with the help of (\ref{A-A}).
Omission points in the fifth line stand for similar evaluations of $\RTr{i}\dots \RTr{2}$,
the resulting expression coincides with (\ref{pair-Ta}).
\hfill$\blacksquare$

\smallskip
For $a_n$ relation (\ref{pair-Ta}) simplifies to
\be
\lb{pair-Tan}
\langle T,\, a_n\rangle\, =\, q^{-1}\id\, ,
\ee
which obviously agrees with (\ref{RE-sl}).
So, we checked a consistency of the normalizations $q^{n-{1\over n}}$
in (\ref{L-iz-univ}) and $\eta=q^{1/n}$ for the
Drinfeld-Jimbo R-matrices with the $SL_q(n)$ reduction condition (\ref{RE-sl}).

\begin{cor}
\lb{cor-A2}
In the conditions of proposition \ref{propA} the pairing $\langle\cdot ,\cdot\rangle$
can be extended for the spectral variables (\ref{Ch-E-hom}):
\be
\lb{pair-Tmu}
\langle T, \, \mu_{\alpha}\rangle\, =\, q^{(2\alpha +2\delta_{\alpha n}-n-{3\over n}-1)}\id , \qquad
\alpha=1,\dots ,n.
\ee
\end{cor}

\ni
{\em Proof.~}
Let us rescale spectral variables $\tilde\mu_\alpha := q^{{3\over n}-2\delta_{\alpha n}}\mu_\alpha$.
For rescaled variables (\ref{pair-Tmu}) reads
\be
\lb{pair-Tmu-tilde}
\langle T, \, \tilde\mu_\alpha\rangle = q^{(2\alpha -n-1)}\id .
\ee
Using q-binomial identity
$\,q^i {n-1 \choose i}_q + q^{i-n}{n-1 \choose i-1}_q= {n\choose i}_q\, $
it is straightforward to derive from (\ref{pair-Tmu-tilde})
pairings of the elementary symmetric functions
$e_i(\tilde\mu)$ by induction on $n$:
$
\langle T, \, e_i(\tilde\mu)\rangle\, =\,{\textstyle {n\choose i}_q}\, .
$
Using (\ref{mu-strih}) it is then straightforward to derive pairings
for elementary symmetric functions in original spectral variables ---
$\langle T, \, e_i(\mu)\rangle$ --- which under identification $a_i\mapsto e_i(\mu)$
coincide with (\ref{pair-Ta}).
\hfill$\blacksquare$

\section{Three lemmas for subsection \ref{subsec3.4}}
\label{append}
Here we collect some technical results which are used
for establishing relation between the spectra of the left and right invariant
vector fields.

\begin{lem}\lb{lem-a1}
{\bf a)} If the {\rm R}-matrix $\R$ is skew invertible
then the following four statements are equivalent:
{\rm ~i)} the matrix $\D$ is invertible;
{\rm ~ii)} the matrix $\C$ is invertible;
{\rm ~iii)}\vspace{-1mm} the {\rm R}-matrix $R^{-1}$ is skew invertible;
{\rm ~iv)} the R-matrix $\Rd:=P\R^{-1} P$ is skew invertible.
One has
\be
\lb{CDR-R}
\Dd\, =\, \Cm\, =\, (\D\,)^{-1}\, , \qquad
\Cd\, =\, \Dm\, =\, (\C\,)^{-1}\, .
\ee

\ni
{\bf b)~} Let $\R$ be the Hecke type R-matrix generating
representations $\Rro$\vspace{-1mm} (\ref{R-rep2}) of the algebras ${\cal H}_k(q)$. Then
the {\rm R}-matrix $\Rd$ is Hecke type as well, and $\Rdro$ are representations of the algebras
${\cal H}_k(q^{-1})$. If additionally the parameter $q$ satisfies conditions {\rm\bf [k]}
(\ref{restrict}) so that  the idempotent
$a^{(k)}|_{q\leftrightarrow q^{-1}}\in {\cal H}_k(q^{-1})$,
 (see (\ref{q-anti1})) is well defined, then
\ba
\lb{RdroA}
\Adro{(k)}& :=& \Rdro(a^{(k)}|_{q\leftrightarrow q^{-1}})\, =\,
\Upsilon^{(k)}_{P}\, A^{(k)}\, \Upsilon^{(k)}_P\, ,
\\[0mm]
\lb{R-Adro}
\Rd_i \Adro{(k)}& =& -q \Adro{(k)}\, , \quad\forall\, i=1,\dots ,k-1\, .
\ea
Here $\Upsilon^{(k)}_P=(\Upsilon^{(k)}_P)^{-1}$ is a particular $R=P$ case of an operator
$\Upsilon^{(k)}_R\in {\rm Aut}(V^{\otimes k})$, defined inductively for any {\rm R}-matrix $\R$
\be
\lb{Theta}
\Upsilon^{(1)}_R\, :=\, 1\, , \quad \Upsilon^{(k+1)}_R\, :=\,
\bigl(\R_1\R_2\dots\R_k\bigr)\,\Upsilon^{(k)}_R\, =\, \Upsilon^{(k)}_R\bigl(\R_k\dots\R_2\R_1\bigr)
\quad \forall\,k=2,3,\dots \, .
\ee
This operator performs reflection of the indices of the {\rm R}-matrices
\be
\lb{R-Theta}
\R_i\, \Upsilon^{(k)}_R\, =\,\Upsilon^{(k)}_R\, \R_{k-i} \qquad\quad\;
\forall\, i,k:\; 1\leq i< k\, .
\ee
The particular element $\Upsilon^{(k)}_P$ enjoys also relations
\ba
\lb{P-Theta}
\R_i\, \Upsilon^{(k)}_P& =&\Upsilon^{(k)}_P\, (P\R P)_{k-i} \quad
\forall\, i,k:\; 1\leq i< k\, ,\quad\forall\;\mbox{{\rm R}-matrix $\R$},
\\
\lb{M-Theta}
M_i\, \Upsilon^{(k)}_P& =&\Upsilon^{(k)}_P\, M_{k-i+1} \qquad
\forall\, i,k:\; 1\leq i\leq k\, ,\quad\forall\,M\in {\rm End}_W(V)\, ,
\ea
where $W$ is an arbitrary $\mathbb C$-linear space.
\end{lem}

\ni
{\em Proof.~}
The first equality in both formulas (\ref{CDR-R}) is
proved in a more general setting
in \cite{OP}, lemma 3.6 c).
The second equality is proved in \cite{I},  section 3.1, proposition 2.
Relations (\ref{R-Theta}) and (\ref{P-Theta}) for matrices $\Upsilon^{(k)}_R$,
$\Upsilon^{(k)}_P$ follow directly from
(\ref{ybe}) and from equalities
\be
\lb{RPP}
\R_1 P_2 P_1\, =\, P_2 P_1 \R_2\, , \quad \R_2 P_1 P_2\, =\, P_1 P_2 \R_1\, .
\ee
Equalities (\ref{M-Theta}) are obvious.
Relations (\ref{R-Adro}) are byproducts of (\ref{RdroA}) and
(\ref{idemp-1}). The second equality in (\ref{RdroA}) follows from
(\ref{q-anti1}), (\ref{q-anti2}), (\ref{P-Theta}), and the Hecke relation for $\Rd$: $\Rd = P\R P-(q-q^{-1})\id$.
\hfill$\blacksquare$

\begin{lem}\lb{lem-a2}
Let $M$ be a matrix of left invariant vector fields for the Hecke type HD algebra
${\goth D\goth G}[\R,\gamma]$.
For any $i\geq 1$ and $j\geq 0$ one has
\ba
\lb{MJ-many}
\lefteqn{\Bigl(\Md{1}\Jd_1\Bigr)\Bigl(\Md{2}\Jd_2\Bigr)\dots \Bigl(\Md{i}\Jd_i\Bigr)
\, I_{i+1,\dots ,i+j}}&&
\\[0mm]
\nn
&&\hspace{-7mm}=
\gamma^{i(i+1)}\, \RTr{i+j+1,\dots ,2i+j}
\left( \Upsilon^{(i+j)}_P\Upsilon^{(2i+j)}_P
(LT)_i\dots (LT)_1 \Upsilon^{(2i)}_{\subRd} (T_i\dots T_1)^{-1}
\Upsilon^{(2i+j)}_P \Upsilon^{(i+j)}_P
\right),
\ea
where $I_{i+1,\dots ,i+j}$ is the identity
operator acting in the component spaces $V$ with labels  $i+1,\dots ,i+j$.
\end{lem}

\ni
{\em Proof.~}
Consider following sequence of transformations
\ba
\nn
\Md{1} \Jd_1 &=& M_1\, =\,
(T^{-1} LT)_1\, =\, \RTr{2} \Bigl(\underline{T^{-1}_1 \R_1 L_1} T_1\Bigr)\, =\,
\gamma^2\, \RTr{2}\Bigl(L_2 \underline{T_1^{-1}\R^{-1} T_1}\Bigr)\,
\\
\lb{tech1}
&=&
\gamma^2\, \RTr{2} \Bigl((LT)_2\, R_1^{-1}\, T_2^{-1}\Bigr)\, =\,
\gamma^2\, \RTr{2} \Bigl(P_1 (LT)_1\, \Rd_1\, T_1^{-1} P_1\Bigr) .
\ea
Here in the first line we transform underlined expressions using (\ref{LT}) and (\ref{RTT}),
and in the last line we apply the definition (\ref{PRP}).
Relation (\ref{tech1}) reproduces formula (\ref{MJ-many}) for $i=1$ and $j=0$.
By a repeated application of formula (c.f. with (\ref{dmat4}))
\be
\lb{PXP}
\RTr{j+1} (P_j X P_j)\, =\, I_j\,\, \RTr{j} (X)
\qquad \forall\, X\in {\rm End}_W(V^{\otimes j}).
\ee
we can rewrite it as (\ref{MJ-many}) with $i=1$ and arbitrary $j>0$
\ba
\lb{tech2}
\Bigl(\Md{1} \Jd_1\Bigr)\, I_{2,\dots j+1}
&=&
\gamma^2\, \RTr{j+2} \Bigl((P_{j+1}\dots P_1) (LT)_1\, \Rd_1\, T_1^{-1} (P_1\dots P_{j+1})\Bigr) .
\ea
In a similar way, for any
value of $i$ relations (\ref{MJ-many}) with $j>0$ follow from
that with $j=0$ by a repeated application of (\ref{PXP}).
Therefore, it is enough considering the case $j=0$.

Using relations (\ref{tech2}) and (\ref{RPP})
we can rewrite an expression $\Md{i}\Jd_i$ in a following way\vspace{-5mm}
\ba
\nn
\lefteqn{\Md{i}\Jd_i\, =\, (\Rd_{i-1}\dots \Rd_1)\, M_{1}\, (\Rd_1\dots \Rd_{i-1})}&&\hspace{6mm}
\\[-2mm]
\lb{MJ}
&= \gamma^2\, \RTr{i+1}\Bigl(
(P_{i}\dots P_2 P_1) (LT)_1\, (\Rd_i\dots \Rd_2\Rd_1\Rd_2\dots \Rd_i)
\, (T_1)^{-1} (P_1 P_2\dots P_{i})\Bigr) .&
\ea

Now we are ready to prove formula (\ref{MJ-many}) by induction on $i$.
Assuming that (\ref{MJ-many}) with $j=1$ is valid for the product of $(i-1)$ factors
we transform the product of $i$ factors
\ba
\nn
&&\hspace{-8mm}\Bigl(\Md{1}\Jd_1\Bigr)\Bigl(\Md{2}\Jd_2\Bigr)\dots \Bigl(\Md{i}\Jd_i\Bigr)=
\gamma^{(i-1)i}\, \RTr{ i+1,\dots 2i-1}\Bigl(\Upsilon^{(i)}_P \Upsilon^{(2i-1)}_P (LT)_{i-1}\dots
(LT)_1 \times
\\
\nn
&&\hspace{71mm}
\times \Upsilon_{\subRd}^{(2i-2)}(T_{i-1}\dots T_1)^{-1} \Upsilon^{(2i-1)}_P \Upsilon^{(i)}_P
\Bigr) \Bigl(\Md{i} \Jd_i\Bigr)
\ea\vspace{-5mm}

\ni
Next, we apply formulas (\ref{P-Theta}), (\ref{M-Theta}) to move the last factor
$(\Md{i}\Jd_i)$\vspace{-2mm}
in this expression leftwards. The result is
\ba
\nn
&&\hspace{-8mm}=
\gamma^{(i-1)i}\, \RTr{ i+1,\dots 2i-1}\Bigl(\Upsilon^{(i)}_P \Upsilon^{(2i-1)}_P (LT)_{i-1}\dots
(LT)_1 \Upsilon_{\subRd}^{(2i-2)}  (T_{i-2}\dots T_1)^{-1}\times
\\
\nn
&&\hspace{85mm}
\times\Bigl(T_1^{-1}(\Md{i}\Jd_i)^{\uparrow 1}\Bigr)^{\uparrow (i-2)}
\Upsilon^{(2i-1)}_P \Upsilon^{(i)}_P \Bigr),
\ea\vspace{-4mm}

\ni
where \vspace{-3mm}we have used identities $(T_{i-1})^{-1}=(T^{-1}_1)^{\uparrow (i-2)}$~
and~ $(\Md{i}\Jd_i)^{\uparrow (i-1)}=((\Md{i}\Jd_i)^{\uparrow 1})^{\uparrow (i-2)}$
to arrange the terms $(T_{i-1})^{-1}$ and $(\Md{i}\Jd_i)^{\uparrow (i-1)}$ in a  suitable
way.\vspace{-3mm} Next, we use formula (\ref{MT-copies}) for their permutation and then,
in a similar way we move term $(\Md{}\Jd)$ to the left of all the terms $(T_\ast)^{-1}\,$:
\ba
\nn
&&\hspace{-5mm}\dots =
\gamma^{(i-1)i+2(i-1)}\, \RTr{ i+1,\dots 2i-1}\Bigl(\Upsilon^{(i)}_P \Upsilon^{(2i-1)}_P (LT)_{i-1}\dots
(LT)_1 \Upsilon_{\subRd}^{(2i-2)}\times
\\
\nn
&&\hspace{73mm}
\times (\Md{2i-1}\Jd_{2i-1}) (T_{i-1}\dots T_1)^{-1}
\Upsilon^{(2i-1)}_P \Upsilon^{(i)}_P \Bigr).
\ea
Now we substitute the expression (\ref{MJ}) for $(\Md{2i-1}\Jd_{2i-1})$\vspace{-2mm}
\ba
\nn
&&\hspace{-5mm} =
\gamma^{i(i+1)}\, \RTr{ i+1,\dots 2i}\Bigl(\Upsilon^{(i)}_P \Upsilon^{(2i-1)}_P (LT)_{i-1}\dots
(LT)_1 \Upsilon_{\subRd}^{(2i-2)}(P_{2i-1}\dots P_1)(LT)_1\times
\\
\nn
&&\hspace{18mm}
\times (\Rd_{2i-1}\dots \Rd_2\Rd_1\Rd_2\dots \Rd_{2i-1})(T_1)^{-1}
(P_1\dots P_{2i-1})
(T_{i-1}\dots T_1)^{-1}
\Upsilon^{(2i-1)}_P \Upsilon^{(i)}_P \Bigr)
\ea
and move the term $(P_{2i-1}\dots P_1)$ leftwards and the term
$(P_1\dots P_{2i-1})$ rightwards close to the terms $\Upsilon_P^{(2i-1)}$.
Finally, using (\ref{Theta})
we complete the calculation
\ba
\nn
&&\hspace{-5mm} =
\gamma^{i(i+1)}\, \RTr{ i+1,\dots 2i}\Bigl(\Upsilon^{(i)}_P \Upsilon^{(2i)}_P (LT)_{i}\dots
(LT)_1 \Upsilon_{\subRd}^{(2i)}(T_{i}\dots T_1)^{-1}
\Upsilon^{(2i)}_P \Upsilon^{(i)}_P \Bigr) , \hspace{23mm}
\ea\vspace{-4mm}

\ni
Here we transformed terms containing $\Rd$ in a following way
\ba
\nn
\Upsilon_{\subRd}^{(2i-2)\uparrow 1}(\Rd_{2i-1}\dots \Rd_2\Rd_1\Rd_2\dots \Rd_{2i-1})\, =\,
\Upsilon_{\subRd}^{(2i-2)\uparrow 1}(\Rd_{1}\dots \Rd_{2i-2}\Rd_{2i-1}\Rd_{2i-2}\dots \Rd_{1})
\hspace{7mm}&&
\\
\nn
=\, (\Rd_{1}\dots \Rd_{2i-1})\Upsilon_{\subRd}^{(2i-2)}(\Rd_{2i-2}\dots \Rd_{1})\, =\,
(\Rd_{1}\dots \Rd_{2i-1})\Upsilon_{\subRd}^{(2i-1)}\,=\,\Upsilon_{\subRd}^{(2i)}.
&&\hspace{0mm}\blacksquare
\ea

\begin{lem}
\lb{lem-a3}
The operators
$J_k$ (\ref{J-k}) and $\Upsilon^{(k)}$ (\ref{Theta}) associated with
a skew invertible {\rm R}-matrix $\R$
satisfy relations
\ba
\lb{Theta-J}
\RTr{i+1,\dots ,2i}\Upsilon_{\R}^{(2i)}\, =\, \Bigl(\Upsilon_\R^{(i)}\Bigr)^4\, =\,
\Bigl(J_1 J_2\dots J_i\Bigr)^2\, .
\ea
\end{lem}

\ni
{\em Proof.~}
Calculation proceeds as follows:
\ba
\nn
\RTr{i+1,\dots ,2i}\Upsilon_{\R}^{(2i)}& =&
\RTr{i+1,\dots ,2i-1}\Bigl(\Upsilon_{\R}^{(2i-1)}\,(\RTr{2i} \R_{2i-1})\,(\R_{2i-2}\dots \R_1) \Bigr)
\\[2mm]
\nn
&=& \RTr{i+1,\dots ,2i-1}\Bigl((\R_1\dots \R_{i-1})\,\Upsilon_{\R}^{(2i-1)}\,(\R_{i-1}\dots \R_1) \Bigr)
\\[2mm]
\nn
\dots &=& (\R_1\dots \R_{i-1})^i\,\Upsilon_{\R}^{(i)}\,(\R_{i-1}\dots \R_1)
(\R_{i-2}\dots \R_1)\dots (\R_2 \R_1) \R_1
\\[2mm]
\nn
&=& \Bigl(J_1 J_2\dots J_i\Bigr)\Bigl(\Upsilon_{\R}^{(i)}\Bigr)^2\, =\, \Bigl(\Upsilon_{\R}^{(i)}\Bigr)^4.
\ea
Here in passing to the second line we calculated the R-trace $\RTr{2i}$\vspace{-1.2mm}
with the help of (\ref{dmat1}) and then used (\ref{R-Theta}) to move
$(i-1)$ R-matrices to the left of the term $\Upsilon_\R^{(2i-1)}$. Expression in the third line
results from a similar calculations of the R-traces $\RTr{2i-1},\; \dots ,\;\RTr{i+1}$,
consecutively. Equalities in the last line result from rearranging factors of
the product $(\R_1\dots \R_{i-1})^i$.
\hfill$\blacksquare$

\newcommand\arxiv[1]{%use this to give preprint # in refs
\href{http://arxiv.org/abs/#1}{\tt arXiv:#1}}

%\addtocontents{toc}{\contentsline{section}{\numberline{}\hspace{-4pt}References}{\pageref{refer}}}

\end{document}